\documentclass{amsart}
\usepackage{xypic}
\usepackage[all]{xy}
\usepackage{url}
\usepackage{amssymb}
\def\half{\frac12}
\let\<\langle
\let\>\rangle
\def\figdir/{/u/levy/MSRI/Book54/main/}

\newcommand\Seqinitial{\xymatrix{
{\Diff^{\infty}(M)}\ar[r]^{P}&{\Diff^{\infty}(M)^k}\ar[r]^{Q}&{\Diff^{\infty}(M)}^r\\
}}
\newcommand\Seqpoint{\xymatrix{
{\Vect^{\infty}(M)}\ar[r]^{DP_{\Id}}&{\Vect^{\infty}(M)^k}\ar[r]^{DQ_{\pi}}&{\Vect^{\infty}(M)^r}\\
}}

\newcommand\inv{^{-1}}

\newcommand\fh{\mathcal H}

\newcommand\ff{\mathfrak F}

\newcommand\G{\Gamma}

\newcommand\fL{\mathfrak L}
\newcommand\g{\gamma}
\newcommand\Ra{\mathbb R}\newcommand\Ta{\mathbb T}
\newcommand\Ca{\mathbb C}
\newcommand\Za{\mathbb Z}
\newcommand\Qa{\mathbb Q}

\newcommand\Ha{\mathbb H}
\newcommand\sw{\mathcal W}
\newcommand\dk{\disp_S}

 \DeclareMathOperator{\Isom}{Isom}
 \DeclareMathOperator{\Id}{Id}
 
\DeclareMathOperator{\Diff}{Diff} 
\DeclareMathOperator{\Aff}{Aff} \DeclareMathOperator{\Vect}{Vect}

\DeclareMathOperator{\Homeo}{Homeo}
 
\DeclareMathOperator{\Hom}{Hom} \DeclareMathOperator{\Aut}{Aut}

\DeclareMathOperator{\disp}{disp}\DeclareMathOperator{\Out}{Out}

\newtheorem{theorem}{Theorem}[section]

\newtheorem{lemma}[theorem]{Lemma}

\newtheorem{defn}[theorem]{Definition}

\newtheorem{question}[theorem]{Question}
\newtheorem{conjecture}[theorem]{Conjecture}

\newtheorem{remarks}{Remarks}
\newtheorem{remark}[remarks]{Remark}

\begin{document}

\title[Local rigidity of group actions]{Local rigidity of
group actions:\\past, present, future}
\author{David Fisher}
\address{Department of Mathematics\\
Rawles Hall\\
Indiana University\\
Bloomington, IN 47405}
\email{fisherdm@indiana.edu}

\thanks{Author partially supported by NSF grant DMS-0226121 and a PSC-CUNY grant.}

\dedicatory{To Anatole Katok on the occasion of his 60th birthday.}

\begin{abstract}
This survey aims to cover the motivation for and history of the
study of local rigidity of group actions.  There is a particularly
detailed discussion of recent results, including outlines of some
proofs. The article ends with a large number of conjectures and
open questions and aims to point to interesting directions for
future research.
\end{abstract}

\maketitle

\section{Prologue}
\label{section:prologue}

Let $\G$ be a finitely generated group, $D$ a topological group,
and $\pi:\G{\rightarrow}D$ a homomorphism.  We wish to study the
space of deformations or perturbations of $\pi$.  Certain trivial
perturbations are always possible as soon as $D$ is not discrete,
namely we can take $d{\pi}d{\inv}$ where $d$ is a small element of
$D$.  This motivates the following definition:

\begin{defn}
\label{definition:locallyrigid} Given a homomorphism
$\pi:\G{\rightarrow}D$, we say $\pi$ is {\em locally rigid} if any
other homomorphism $\pi'$ which is close to $\pi$ is conjugate to
$\pi$ by a small element of $D$.
\end{defn}

\noindent We topologize $\Hom(\G,D)$ with the compact open
topology which means that two homomorphisms are close if and only
if they are close on a generating set for $\G$.  If $D$ is path
connected, then we can define {\em deformation rigidity} instead,
meaning that any continuous path of representations $\pi_t$
starting at $\pi$ is conjugate to the trivial path $\pi_t=\pi$ by
a continuous path $d_t$ in $D$ with $d_0$ being the identity in
$D$.  If $D$ is an algebraic group over $\Ra$ or $\Ca$, it is
possible to prove that deformation rigidity and local rigidity are
equivalent since $\Hom(\G,D)$ is an algebraic variety and the
action of $D$ by conjugation is algebraic; see \cite{Mu}, for example. For
$D$ infinite dimensional and path-connected, this equivalence is
no longer clear.

The study of local rigidity of lattices in semi-simple Lie groups
is probably the beginning of the general study of rigidity in
geometry and dynamics, a subject that is by now far too large for
a single survey. See \cite{Sp1} for the last attempt at a
comprehensive survey and \cite{Sp2} for a more narrowly focused
updating of that survey.  Here we abuse language slightly by
saying a subgroup is locally rigid if the defining embedding is
locally rigid as a homomorphism.  See subsection
\ref{subsection:selbergweil} for a brief history of local rigidity
of lattices and some discussion of subsequent developments that
are of particular interest in the study of rigidity of group
actions.

In this article we will focus on a survey of local rigidity when
$D=\Diff^{\infty}(M)$ or occasionally $D=\Diff^k(M)$ for some
finite $k$.  Here we often refer to a homomorphism
$\pi:\G{\rightarrow}\Diff^{\infty}(M)$ as an action, since it can
clearly be thought of as $C^{\infty}$ action
$\pi:\G{\times}M{\rightarrow}M$.  We will consistently use the
same letter $\pi$ to denote either the action or the homomorphism
to $\Diff^{\infty}(M)$. The title of this article refers to this
interpretation of $\pi$ as defining a group action. In this
context, one considers rigidity of actions of connected groups as
well as of discrete groups. In cases where $\G$ has any topology,
we will always only study continuous actions, i.e. ones for which
the defining homomorphism $\pi$ is a continuous map.

One can, in this context, develop more refined notions of local
rigidity, since the topology on $\Diff^{\infty}(M)$ is an inverse
limit topology. This means that two $C^{\infty}$ diffeomorphisms
of $M$ are close precisely when they are $C^k$ close for some
large value of $k$.  The most exhaustive definition of local
rigidity is probably the following:

\begin{defn}
\label{defn:locallyrigidrefined} Let $\G$ be a discrete group and
$\pi:{\G}{\rightarrow}\Diff^{k}(M)$ a homomorphism where $k$ is
either a positive integer or $\infty$. We say that $\pi$ is {\em
$C^{k,l,i,j,m}$ locally rigid} if any
$\pi':\G{\rightarrow}\Diff^{l}(M)$ which is close to $\pi$ in the
$C^i$ topology is conjugate to $\pi$ by a $C^j$ diffeomorphism
$\phi$ which is $C^m$ small.  Here $l,i,j,m$ are all either
non-negative integers or $\infty$ and the only a priori constraint
are $i{\leq}\min(k,l)$ and $m{\leq}j$.  When $j=0$, we will call
the action {\em stable} or {\em structurally stable}.  When $j>0$,
we will call the action {\em locally rigid} or simply {\em rigid}.
\end{defn}

\noindent We will avoid using this cumbersome notation when at all
possible. There is a classical, dynamical notion of structural
stability which is equivalent to $C^{1,1,1,0,0}$ local rigidity.
I.e. a $C^1$ action $\pi$ of a group $\G$ is {\em structurally
stable} if any $C^1$ close $C^1$ action of $\G$ is conjugate to
$\pi$ by a small homeomorphism. For actions of $\Za$ this notion
arose in hyperbolic dynamics in the work of Anosov and Smale
\cite{An,Sm}. From a dynamical point of view structural stability
is important since it allows one to control dynamical properties
of an open set of actions in $\Diff^1(M)$.  Local rigidity can be
viewed as a strengthening of this property in that it shows that
an open set of actions is exhausted by smooth conjugates of a
single action.

Though actions of $\Za$ and free groups on $k$ generators are
often structurally stable, they are never locally rigid, and it is
an interesting question as to how ``large" a group needs to be in
order to have actions which are locally rigid. Many of the
original questions and theorems concerning local rigidity were for
lattices in higher rank semi-simple Lie groups, where here higher
rank means that all simple factors have real rank at least $2$.
(See subsection \ref{subsection:groups} for a definition of rank.)
Fairly early in the theory it became clear that local rigidity
often held, and was in fact easier to prove, for certain actions
of higher rank abelian groups, i.e. $\Za^k$ for $k{\geq}2$, see
\cite{KL1}. In addition, local rigidity results have been proven
for actions of a wider variety of groups, including
\begin{enumerate} \item certain
non-volume preserving actions of lattices in $SO(1,n)$ in
\cite{Kan1} \item all isometric actions of groups with property
$(T)$ in \cite{FM2}, \item certain affine actions of lattices in
$SP(1,n)$ in \cite{Hi}.
\end{enumerate}
\noindent There is extremely interesting related work of Ghys,
older than the work just mentioned, which shows that the space of
deformations of certain actions of surface groups on the circle is
finite dimensional \cite{Gh1,Gh2,Gh3}. Ghys also proved some very
early results on local and global rigidity of very particular
actions of connected solvable groups, see \cite{GhS,Gh1} and
subsection \ref{subsection:boundaries}.

The study of local rigidity of group actions has had three primary
historical motivations, one from the theory of lattices in Lie
groups, one from dynamical systems theory and a third from the
theory of foliations. (This statement is perhaps a bit coarse, and
there is heavy overlap between these motivations, particularly the
second and the third.) The first is the general study of rigidity
of actions of large groups, as discussed in \cite{Z3,Z4}, see
\cite{La,FK} for more up to date surveys. This area is motivated
by the study of rigidity of lattices in semi-simple Lie groups,
particularly by Margulis' super-rigidity theorem and it's
non-linear generalization by Zimmer to a cocycle super-rigidity
theorem, see subsection \ref{subsection:selbergweil} and \cite{Z4}
for more discussion. This motivation also stems from an analogy
between semi-simple Lie groups and diffeomorphism groups.  When
$M$ is a compact manifold, not only is $\Diff^{\infty}(M)$ an
infinite dimensional Lie group, but its connected component is
simple. Simplicity of the connected component of
$\Diff^{\infty}(M)$ was proven by Thurston using results of
Epstein and Herman \cite{Th2,Ep,Hr}. Herman had used Epstein's
work to see that the connected component of
$\Diff^{\infty}(\Ta^n)$ is simple and Thurston's proof of the
general case uses this. See also further work on the topic by
Banyaga and Mather \cite{Ba1,Mt1,Mt2,Mt3}, as well as Banyaga's
book \cite{Ba2}.

The dynamical motivation for studying rigidity of group actions
comes from the study of structural stability of diffeomorphisms
and flows in hyperbolic dynamics, see the introduction of
\cite{KS2}. This area provides many of the basic techniques by
which results in the area have been proven, particularly early in
the history of the field. Philosophically, hyperbolic
diffeomorphisms are structurally stable, group actions generated
by structurally stable diffeomorphisms are quite often
structurally stable, and the presence of a large group action
frequently allows one to improve the regularity of the conjugacy.
See subsection \ref{subsection:hyperbolic} for a brief history of
relevant results on structural stability and subsections
\ref{subsection:volumepreserving}, \ref{subsection:boundaries},
and \ref{subsection:fmsketch} for some applications of these
ideas.

The third motivation for studying rigidity of group actions comes
from the theory of foliations.  Many techniques and ideas in this
area are also related to work on hyperbolic dynamics, and many of
the foliations of interest are dynamical foliations of hyperbolic
dynamical systems. A primary impetus in this area is the theory of
codimension one foliations, and so many of the ideas here were
first developed either for groups acting on the circle or for
actions of connected groups on manifolds only one dimension larger
then the acting group.  See particularly \cite{GhS,Gh1} for the
early history of these developments.

\subsection*{Some remarks on biases and omissions.} Like any
survey of this kind, this work is informed by it's authors biases
and experiences.  The most obvious of these is that my point of
view is primarily motivated by the study of rigidity properties of
semi-simple Lie groups and their lattices, rather than primarily
motivated by hyperbolic dynamics or foliation theory.  This
informs the biases of this article and a very different article
would result from different biases.

There are two large omissions in this article.  The first omission
is that it is primarily occupied with local rigidity of discrete
group actions.  When similar results are known for actions of Lie
groups, they are mentioned, though frequently only special cases
are stated. This is partially because results in this context are
often complicated by the need to consider time changes, and I did
not want to dwell on that issue here.  The second omission is that
little to no care is taken to state optimal results relating the
various constants in $C^{k,l,i,j,m}$ local rigidity. Dwelling on
issues of regularity seemed likely to obscure the main line of the
developments, so many results are stated without any explicit
mention of regularity.  Usually this is done only when the action
can be shown to be {\em locally rigid} in $\Diff^{\infty}(M)$ in
the sense of Definition \ref{definition:locallyrigid}. This
implicitly omits both the degree of regularity to which the
perturbed and unperturbed actions are close and the degree of
regularity with which the size of the conjugacy is small. In other
words local rigidity is $C^{\infty,\infty,i,\infty,m}$ local
rigidity for some unspecified $i$ and $m$, and I always fail to
specify $i$ and $m$ even when they are known. Occasionally a
result is stated that only produces a finite regularity conjugacy,
with this issue only remarked on following the statement of the
result. It seems quite likely that most results of this kind can
be improved to produce $C^{\infty}$ conjugacies using the
techniques of \cite{FM2,FM3}, see discussion at the end of
Section \ref{subsection:fmsketch}.

Lastly we remark that the study of local rigidity of group actions
is often closely intertwined with the study of {\em global
rigidity} of group actions.  The meaning of the phrase global
rigidity is not entirely precise, but it is typically used to
cover settings in which one can classify all group actions
satisfying certain hypotheses on a certain manifold or class of
manifolds.  The study of global rigidity is too broad and
interesting to summarize briefly, but some examples are mentioned
below when they are closely related to work on local rigidity. See
\cite{FK,La} for recent surveys concerning both local and global
rigidity.

\section{A brief digression: some examples of groups and actions}
\label{section:examples}

In this section we briefly describe some of the groups that will
play important roles in the results discussed here.  The reader
already familiar with semi-simple Lie groups and their lattices
may want to skip to the second subsection where we give
descriptions of group actions.

\subsection{Semi-simple groups and their lattices.}
\label{subsection:groups}

By a simple Lie group, we mean a connected Lie group all of whose
normal subgroups are discrete, though we make the additional
convention that $\Ra$ and $S^1$ are not simple. By a semi-simple
Lie group we mean the quotient of a product of simple Lie groups
by some subgroup of the product of their centers.  Note that with
our conventions, the center of a simple Lie group is discrete and
is in fact the maximal normal subgroup. There is an elaborate
structure theory of semi-simple Lie groups and the groups are
competely classified, see \cite{He} or \cite{Kn} for details. Here
we merely describe some examples, all of which are matrix groups.
All connected semisimple Lie groups are discrete central
extensions of matrix groups, so the reader will lose very little
by always thinking of matrix groups.
\begin{enumerate} \item The groups $SL(n,\Ra), SL(n,\Ca)$ and
$SL(n,\Ha)$ of $n$ by $n$ matrices of determinant one over the
real numbers, the complex numbers or the quaternions. \item The
group $SP(2n,\Ra)$ of $2n$ by $2n$ matrices of determinant one
which preserve a real symplectic form on $\Ra^{2n}$.  \item The
groups $SO(p,q), SU(p,q)$ and $SP(p,q)$ of matrices which preserve
inner products of signature $(p,q)$ where the inner product is
real linear on $\Ra^{p+q}$, hermitian on $\Ca^{p+q}$ or
quaternionic hermitian on $\Ha^{p+q}$ respectively.
\end{enumerate}

Let $G$ be a semi-simple Lie group which is a subgroup of
$GL(n,\Ra)$.  We say that $G$ has {\em real rank} $k$ if $G$ has a
$k$ dimensional abelian subgroup which is conjugate to a subgroup
of the real diagonal matrices and no $k+1$ dimensional abelian
subgroups with the same property.  The groups in $(1)$ have rank
$n-1$, the groups in $(2)$ have rank $n$ and the groups in $(3)$
have rank $\min(p,q)$.

Since this article focuses primarily on finitely generated groups,
we are more interested in discrete subgroups of Lie groups than in
the Lie groups themselves.  A discrete subgroup $\G$ in a Lie
group $G$ is called a lattice if $G/{\G}$ has finite Haar measure.
The lattice is called {\em cocompact} or {\em uniform} if $G/{\G}$
is compact and {\em non-uniform} or simply not cocompact
otherwise.  If $G=G_1{\times}{\cdots}{\times}G_n$ is a product
then we say a lattice $\G<G$ is {\em irreducible} if it's
projection to each $G_i$ is dense.  More generally we make the
same definition for an {\em almost direct product}, by which we
mean a direct product $G$ modulo some subgroup of the center
$Z(G)$. Lattices in semi-simple Lie groups can always be
constructed by arithmetic methods, see \cite{Bo} and also
\cite{Mr} for more discussion.  In fact, one of the most important
results in the theory of semi-simple Lie groups is that if $G$ is
a semi-simple Lie group without compact factors, then all
irreducible lattices in $G$ are arithmetic unless $G$ is locally
isomorphic to $SO(1,n)$ or $SU(1,n)$.  For $G$ of real rank at
least two, this is Margulis' arithmeticity theorem, which he
deduced from his super-rigidity theorems \cite{Ma1,Ma2,Ma3}. For
non-uniform lattices, Margulis had an earlier proof which does not
use the superrigidity theorems, see \cite{Ma0,Ma00}. This earlier
proof depends on the study of dynamics of unipotent elements on
the space $G/{\Gamma}$, and particularly on what is now known as
the ``non-divergence of unipotent flows". Special cases of the
super-rigidity theorems were then proven for $Sp(1,n)$ and
$F_4^{-20}$ by Corlette and Gromov-Schoen, which sufficed to imply
the statement on arithmeticity given above \cite{Co,GS}. As we
will be almost exclusively concerned with arithmetic lattices, we
do not give examples of non-arithmetic lattices here, but refer
the reader to \cite{Ma3} and \cite{Mr} for more discussion.  A
formal definition of arithmeticity, at least when $G$ is algebraic
is:

\begin{defn}
\label{defn:arithmetic} Let $G$ be a semisimple algebraic Lie
group and $\G<G$ a lattice.  Then $\G$ is arithmetic if there
exists a semi-simple algebraic Lie group $H$ defined over $\Qa$
such that
\begin{enumerate}
\item there is a homomorphism $\pi:H^0{\rightarrow}G$ with compact
kernel, \item there is a rational structure on $H$ such that the
projection of the integer points of $H$ to $G$ are commensurable
to $\G$, i.e. $\pi(H(\Za)){\cap}\G$ is of finite index in both
$H(\Za)$ and $\G$.
\end{enumerate}
\end{defn}

We now give some examples of arithmetic lattices.  The simplest is
to take the integer points in a simple (or semi-simple) group $G$
which is a matrix group, e.g. $SL(n,\Za)$ or $Sp(n,\Za)$. This
exact construction always yields lattices, but also always yields
non-uniform lattices.  In fact the lattices one can construct in
this way have very special properties because they will contain
many unipotent matrices.  If a lattice is cocompact, it will
necessarily contain no unipotent matrices.  The standard trick for
understanding the structure of lattices in $G$ which become
integral points after passing to a compact extension is called
{\em change of base}.  For much more discussion see
\cite{Ma3,Mr,Z2}.  We give one example to illustrate the process.
Let $G=SO(m,n)$ which we view as the set of matrices in
$SL(n+m,\Ra)$ which preserve the inner product
$$\<v,w\>=\biggl(-\sqrt{2}\sum_{i=1}^mv_iw_i\biggr)
+\biggl(\sum_{i=m+1}^{n+m}v_iw_i\biggr)$$
where $v_i$ and $w_i$ are the $i$th components of $v$ and $w$.
This form, and therefore $G$, are defined over the field
$\Qa(\sqrt{2})$ which has a Galois conjugation $\sigma$ defined by
$\sigma(\sqrt{2})=-\sqrt{2}$.  If we looks at the points
$\G=G(\Za[\sqrt{2}])$, we can define an embedding of $\G$ in
$SO(m,n){\times}SO(m+n)$ by taking $\g$ to $(\g,\sigma(\g))$. It
is straightforward to check that this embedding is discrete.  In
fact, this embeds $\G$ in $H=SO(m,n){\times}SO(m+n)$ as integral
points for the rational structure on $H$ where the rational points
are exactly the points $(m,\sigma(m))$ where
$m{\in}G(\Qa(\sqrt{2}))$.  This makes $\G$ a lattice in $H$ and it
is easy to see that $\G$ projects to a lattice in $G$, since $G$
is cocompact in $H$.  What is somewhat harder to verify is that
$\G$ is cocompact in $H$, for which we refer the reader to the
list of references above.

Similar constructions are possible with $SU(m,n)$ or $SP(m,n)$ in
place of $SO(m,n)$ and also with more simple factors and fields
with more Galois automorphisms.  There are also a number of other
constructions of arithmetic lattices using division algebras. See
\cite{Mr} for a comprehensive treatment.

We end this section by defining a key property of many semisimple
groups and their lattices.  This is property $(T)$ of Kazhdan, and
was introduced by Kazhdan in \cite{Ka1} in order to prove that
non-uniform lattices in higher rank semi-simple Lie groups are
finitely generated and have finite abelianization. It has played a
fundamental role in many subsequent developments.  We do not give
Kazhdan's original definition, but one which was shown to be
equivalent by work of Delorme and Guichardet \cite{De,Gu}.

\begin{defn}
\label{definition:T} A group $\G$ has property $(T)$ of Kazhdan if
$H^1(\G,\pi)=0$ for every continuous unitary representation $\pi$
of $\G$ on a Hilbert space. This is equivalent to saying that any
continuous isometric action of $\G$ on a Hilbert space has a fixed
point.
\end{defn}

\begin{remarks}
\begin{enumerate} \item Kazhdan's
definition is that the trivial representation is isolated in the
unitary dual of $\G$.

\item If a continuous group $G$ has property $(T)$ so does any
lattice in $G$.  This result was proved in \cite{Ka1}.

\item Any semi-simple Lie group has property $(T)$ if and only if
it has no simple factors locally isomorphic to $SO(1,n)$ or
$SU(1,n)$. For a discussion of this fact and attributions, see
\cite{HV}.  For groups with all simple factors of real rank at
least three, this is proven in \cite{Ka1}.

\item No noncompact amenable group, and in particular no noncompact
abelian group, has property $(T)$.  An easy averaging argument
shows that all compact groups have property $(T)$.
\end{enumerate}
\end{remarks}

\noindent Groups with property $(T)$ play an important role in
many areas of mathematics and computer science.

\subsection{Some actions of groups and lattices.}
\label{subsection:actions}

Here we define and give examples of the general classes of actions
for which local rigidity results have been proven.  Let $H$ be a
Lie group and $L<H$ a closed subgroup.  Then a diffeomorphism $f$
of $H/L$ is called {\em affine} if there is a diffeomorphism
$\tilde f$ of $H$ such that $f([h])=\tilde f(h)$ where $\tilde f =
A{\circ}\tau_h$ with $A$ an automorphism of $H$ with $A(L)=L$ and
$\tau_h$ is left translation by some $h$ in $H$.  Two obvious
classes of affine diffeomorphisms are left translations on any
homogeneous space and either linear automorphisms of tori or more
generally automorphisms of nilmanifolds.  A group action is called
{\em affine} if every element of the group acts by an affine
diffeomorphism. It is easy to check that the full group of affine
diffeomorphisms $\Aff(H/L)$ is a finite dimensional Lie group and
an affine action of a group $D$ is a homomorphism
$\pi:D{\rightarrow}\Aff(H/L)$. The structure of $\Aff(H/L)$ is
surprisingly complicated in general, it is a quotient of a
subgroup of the group $\Aut(H){\ltimes}H$ where $\Aut(H)$ is a the
group of automorphisms of $H$.  For a more detailed discussion of
this relationship, see \cite[Section 6]{FM1}.  While it is not
always the case that any affine action of a group $D$ on $H/L$ can
be described by a homomorphism
$\pi:D{\rightarrow}\Aut(H){\ltimes}H$, this is true for two
important special cases:
\begin{enumerate}
\item $D$ is a connected semi-simple Lie group and $L$ is a
cocompact lattice in $H$, \item $D$ is a lattice in a semi-simple
Lie group $G$ where $G$ has no compact factors and no simple
factors locally isomorphic to $SO(1,n)$ or $SU(1,n)$, and $L$ is a
cocompact lattice in $H$.
\end{enumerate}
\noindent These facts are \cite[Theorem 6.4 and 6.5]{FM1} where
affine actions as in $(1)$ and $(2)$ above are classified.

The most obvious examples of affine actions of large groups are of
the following forms, which are frequently referred to as {\em
standard actions}:
\begin{enumerate}
\item Actions of groups by automorphisms of nilmanifolds. I.e. let
$N$ be a simply connected nilpotent group, $\Lambda<N$ a lattice
(which is necessarily cocompact) and assume a finitely generated
group $\G$ acts by automorphisms of $N$ preserving $\Lambda$.  The
most obvious examples of this are when $N=\Ra^n$, $\Lambda=\Za^n$
and $\G<SL(n,\Za)$, in which case we have a linear action of $\G$
on $\Ta^n$.

\item Actions by left translations.  I.e. let $H$ be a Lie group
and $\Lambda<H$ a cocompact lattice and $\G<H$ some subgroup. Then
$\G$ acts on $H/\Lambda$ by left translations.  Note that in this
case $\G$ need not be discrete.

\item Actions by isometries.  Here $K$ is a compact group which
acts by isometries on some compact manifold $M$ and $\G<K$ is a
subgroup.  Note that here $\G$ is either discrete or a discrete
extension of a compact group.
\end{enumerate}

We now briefly define a few more general classes of actions, for
which local rigidity results are either known or conjectured. We
first fix some notations. Let $A$ and $D$ be topological groups,
and $B<A$ a closed subgroup. Let
$\rho:D{\times}A/B{\rightarrow}A/B$ be a continuous affine action.

\begin{defn}
\label{definition:affine}
\begin{enumerate}
\item
Let $A,B,D$ and $\rho$
be as above. Let $C$ be a compact group of affine diffeomorphisms
of $A/B$ that commute with the $D$ action. We call the action of
$D$ on $C{\backslash}A/B$ a {\em generalized affine action}.
\item
Let $A$, $B$, $D$ and $\rho$ be as in $1$ above.
Let $M$ be a compact Riemannian manifold and
$\iota:D{\times}A/B{\rightarrow}\Isom(M)$ a $C^1$ cocycle.  We
call the resulting skew product $D$ action on $A/B{\times}M$ a
{\em quasi-affine action}. If $C$ and $D$ are as in $2$, and
$\alpha:D{\times}C{\backslash}A/B{\rightarrow}\Isom(M)$ is a $C^1$
cocycle, then we call the resulting skew product $D$ action on
$C{\backslash}A/B{\times}M$ a {\em generalized quasi-affine
action}.
\end{enumerate}
\end{defn}

For many of the actions we consider here, there will be a
foliation of particular importance.  If $\rho$ is an action of a
group $D$ on a manifold $N$, and $\rho$ preserves a foliation
$\ff$ and a Riemannian metric along the leaves of $\ff$, we call
$\ff$ a {\em central foliation} for $\rho$.  For quasi-affine and
generalized quasi-affine actions on manifolds of the form
$C{\backslash}A/B{\times}M$ the foliation by leaves of the
$\{[a]\}{\times}M$ is always a central foliation.  There are also
actions with more complicated central foliations.  For example if
$H$ is a Lie group, $\Lambda<H$ is discrete and a subgroup $G<H$
acts on $H/{\Lambda}$ by left translations, then the foliation of
$H/{\Lambda}$ by orbits of the centralizer $Z_H(G)$ of $G$ in $H$
is a central foliation. It is relatively easy to construct
examples where this foliation has dense leaves.  Another example
of an action which has a foliation with dense leaves is to embed
the $\Za[\sqrt{2}]$ points of $SO(m,n)$ into $SL(2(m+n),\Za)$ as
described in the preceding subsection and then let this group act
on $\Ta^{2(m+n)}$ linearly. It is easy to see in this case that
the maximal central foliation for the action is a foliation by
densely embedded leaves none of which are compact.

We end this section by describing briefly the standard
construction of an {\em induced or suspended action}.  This notion
can be seen as a generalization of the construction of a flow
under a function or as an analogue of the more algebraic notion of
inducing a representation.  Given a group $H$, a (usually closed)
subgroup $L$, and an action $\rho$ of $L$ on a space $X$, we can
form the space $(H{\times}X)/L$ where $L$ acts on $H{\times}X$ by
$h\cdot(l,x)=(lh{\inv},\rho(h)x)$.  This space now has a natural
$H$ action by left multiplication on the first coordinate.  Many
properties of the $L$ action on $X$ can be studied more easily in
terms of properties of the $H$ action on $(H{\times}X)/L$.  This
construction is particularly useful when $L$ is a lattice in $H$.

\section{Pre-history}
\label{section:avantlalettre}

\subsection{Local and global rigidity of homomorphisms into finite dimensional groups.}\label{subsection:selbergweil}

The earliest work on local rigidity in the context of Definition
\ref{definition:locallyrigid} was contained in series of works by
Calabi--Vesentini, Selberg, Calabi and Weil, which resulted in the
following:

\begin{theorem}
\label{theorem:rigidityoflattices} Let $G$ be a semi-simple Lie
group and assume that $G$ is not locally isomorphic to
$SL(2,\Ra)$. Let $\G<G$ be an irreducible cocompact lattice, then
the defining embedding of $\G$ in $G$ is locally rigid.
\end{theorem}

\begin{remarks}
\begin{enumerate}   \item If $G=SL(2,\Ra)$
the theorem is false and there is a large, well studied space of
deformation of $\G$ in $G$, known as the Teichmuller space. \item
There is an analogue of this theorem for lattices that are not
cocompact.  This result was proven later and has a more
complicated history which we omit here. In this case it is also
necessary to exclude $G$ locally isomorphic to $SL(2,\Ca)$.
\end{enumerate}
\end{remarks}

This theorem was originally proven in special cases by Calabi,
Calabi--Vesentini and Selberg.   In particular, Selberg gives a
proof for cocompact lattices in $SL(n,\Ra)$ for $n\geq{3}$ in
\cite{S}, Calabi--Vesentini give a proof when the associated
symmetric space $X=G/K$ is K\"{a}hler in \cite{CV} and Calabi
gives a proof for $G=SO(1,n)$ where $n{\geq}3$ in \cite{C}.
Shortly afterwards, Weil gave a complete proof of Theorem
\ref{theorem:rigidityoflattices} in \cite{We1,We2}.

In all of the original proofs, the first step was to show that any
perturbation of $\G$ was discrete and therefore a cocompact
lattice. This is shown in special cases in \cite{C,CV,S} and
proven in a somewhat broader context than Theorem
\ref{theorem:rigidityoflattices} in \cite{W1}.

The different proofs of cases of Theorem
\ref{theorem:rigidityoflattices} are also interesting in that
there are two fundamentally different sets of techniques employed
and this dichotomy continues to play a role in the history of
rigidity.   Selberg's proof essentially combines algebraic facts
with a study of the dynamics of iterates of matrices.  He makes
systematic use of the existence of singular directions, or Weyl
chamber walls, in maximal diagonalizable subgroups of $SL(n,\Ra)$.
Exploiting these singular directions is essential to much later
work on rigidity, both of lattices in higher rank groups and of
actions of abelian groups.   It seems possible to generalize
Selberg's proof to the case of $G$ an $\Ra$-split semi-simple Lie
group with rank at least $2$. Selberg's proof, which depended on
asymptotics at infinity of iterates of matrices, inspired Mostow's
explicit use of boundaries in his proof of strong rigidity
\cite{Mo2}. Mostow's work in turn provided inspiration for the use
of boundaries in later work of Margulis, Zimmer and others on
rigidity properties of higher rank groups.

The proofs of Calabi, Calabi--Vesentini and Weil involve studying
variations of geometric structures on the associated locally
symmetric space. The techniques are analytic and use a variational
argument to show that all variations of the geometric structure
are trivial.  This work is a precursor to much work in geometric
analysis studying variations of geometric structures and also
informs later work on proving rigidity/vanishing of harmonic forms
and maps. The dichotomy between approaches based on
algebra/dynamics and approaches that are in the spirit of
geometric analysis  continues through much of the history of
rigidity and the history of local rigidity of group actions in
particular.

Shortly after completing this work, Weil discovered a new
criterion for local rigidity \cite{We3}.  In the context of
Theorem \ref{theorem:rigidityoflattices}, this allows one to avoid
the step of showing that a perturbation of $\G$ remains discrete.
In addition, this result opened the way for understanding local
rigidity of more general representations of discrete groups.

\begin{theorem}
\label{theorem:weil} Let $\G$ be a finitely generated group, $G$ a
Lie group and $\pi:\G{\rightarrow}G$ a homomorphism.  Then $\pi$
is locally rigid if $H^1(\G,\mathfrak{g})=0$.  Here $\mathfrak{g}$
is the Lie algebra of $G$ and $\G$ acts on $\mathfrak{g}$ by
$Ad_G{\circ}\pi$.
\end{theorem}

\noindent  Weil's proof of this result uses only the implicit
function theorem and elementary properties of the Lie group
exponential map.  The same theorem is true if $G$ is an algebraic
group over any local field of characteristic zero.  In \cite{We3},
Weil remarks that if $\G<G$ is a cocompact lattice and $G$
satisfies the hypothesis of Theorem
\ref{theorem:rigidityoflattices}, then the vanishing of
$H^1(\G,\mathfrak{g})$ can be deduced from the computations in
\cite{We2}.  The vanishing of $H^1(\G,\mathfrak{g})$ is proven
explicitly by Matsushima and Murakami in \cite{MM}.

Motivated by Weil's work and other work of Matsushima, conditions
for vanishing of $H^1(\G,\mathfrak{g})$ were then studied by many
authors.  See particularly \cite{MM} and \cite{Rg1}.  The results
in these papers imply local rigidity of many linear
representations of lattices.

To close this section, I will briefly discuss some subsequent
developments concerning rigidity of lattices in Lie groups that
motivated the study of both local and global rigidity of group
actions.

The first remarkable result in this direction is Mostow's rigidity
theorem, see \cite{Mo1} and references there.  Given $G$ as in
Theorem \ref{theorem:rigidityoflattices}, and two irreducible
cocompact lattices $\G_1$ and $\G_2$ in $G$, Mostow proves that
any isomorphism from $\G_1$ to $\G_2$ extends to an isomorphism of
$G$ with itself.  Combined with the principal theorem of
\cite{We1} which shows that a perturbation of a lattice is again a
lattice, this gives a remarkable and different proof of Theorem
\ref{theorem:rigidityoflattices}, and Mostow was motivated by the
desire for a ``more geometric understanding" of Theorem
\ref{theorem:rigidityoflattices} \cite{Mo1}. Mostow's theorem is
in fact a good deal stronger, and controls not only homomorphisms
$\G{\rightarrow}G$ near the defining homomorphism, but any
homomorphism into any other simple Lie group $G'$ where the image
is lattice. As mentioned above, Mostow's approach was partially
inspired by Selberg's proof of certain cases of Theorem
\ref{theorem:rigidityoflattices}, \cite{Mo2}.  A key step in
Mostow's proof is the construction of a continuous map between the
geometric boundaries of the symmetric spaces associated to $G$ and
$G'$.  Boundary maps continue to play a key role in many
developments in rigidity theory. A new proof of Mostow rigidity,
at least for $G_i$ of real rank one, was provided by Besson,
Courtois and Gallot.  Their approach is quite different and has
had many other applications concerning rigidity in geometry and
dynamics, see e.g. \cite{BCG, CF}.

The next remarkable result in this direction is Margulis'
superrigidity theorem. Margulis proved this theorem as a tool to
prove arithmeticity of irredudicible uniform lattices in groups of
real rank at least $2$.  For irreducible lattices in semi-simple
Lie groups of real rank at least $2$, the superrigidity theorems
classifies all finite dimensional linear representations.
Margulis' theorem holds for irreducible lattices in semi-simple
Lie groups of real rank at least two. Given a lattice $\G<G$ where
$G$ is simply connected, one precise statement of some of Margulis
results is to say that any linear representation $\sigma$ of $\G$
{\em almost extends} to a linear representation of $G$.  By this
we mean that there is a linear representation $\tilde \sigma$ of
$G$ and a bounded image representation $\bar \sigma$ of $\G$ such
that $\sigma(\g)=\tilde \sigma(\g)\bar \sigma(\g)$ for all $\g$ in
$G$. Margulis' theorems also give an essentially complete
description of the representations $\bar \sigma$, up to some
issues concerning finite image representations.  The proof here is
partially inspired by Mostow's work: a key step is the
construction of a measurable ``boundary map".  However the methods
for producing the boundary map in this case are very dynamical.
Margulis' original proof used Oseledec Multiplicative Ergodic
Theorem.  Later proofs were given by both Furstenberg and Margulis
using the theory of group boundaries as developed by Furstenberg
from his study of random walks on groups \cite{Fu1,Fu2}.
Furstenberg's probabalistic version of boundary theory has had a
profound influence on many subsequent developments in rigidity
theory. For more discussion of Margulis' superrigidity theorem,
see \cite{Ma1,Ma2,Ma3}.

A main impetus for studying rigidity of group actions on manifolds
came from Zimmer's theorem on superrigidity for cocycles.  This
theorem and it's proof were strongly motivated by Margulis' work.
In fact, Margulis' theorem is Zimmer's theorem for a certain
coccyle $\alpha:G{\times}G/{\Gamma}{\rightarrow}\G$ and the proof
of Zimmer's theorem is quite similar to the proof of Margulis'. In
order to avoid technicalities, we describe only a special case of
this result. Let $M$ be a compact manifold, $H$ a matrix group and
$P=M{\times}H$. Now let a group $\G$ act on $M$ and $P$
continuously, so that the projection from $P$ to $M$ is
equivariant. Further assume that the action on $M$ is measure
preserving and ergodic. If $\G$ is a lattice in a simply
connected, semi-simple Lie group $G$ all of whose simple factors
have real rank at least two then there is a measurable map
$s:M{\rightarrow}H$, a representation $\pi:G{\rightarrow}H$, a
compact subgroup $K<H$ which commutes with $\pi(G)$ and a
measurable map $\G{\times}M{\rightarrow}K$ such that
\begin{equation}
\label{equation:csr} \g{\cdot}s(m)=k(m,\g)\pi(\g)s(\g{\cdot}m).
\end{equation}
\noindent It is easy to check from this equation that the map $K$
satisfies a certain equation that makes it into a cocycle over the
action of $\G$.  One should view $s$ as providing coordinates on
$P$ in which the $\G$ action is {\em almost a product}. For more
discussion of this theorem the reader should see any of
\cite{Fe1,Fe2,FM1,Fu,Z2}. (The version stated here is only proven
in \cite{FM1}, previous proofs all yielded somewhat more
complicated statements.) As a sample application, let $M=\Ta^n$
and let $P$ be the frame bundle of $M$, i.e. the space of frames
in the tangent bundle of $M$. Since $\Ta^n$ is parallelizable, we
have $P=\Ta^n{\times}GL(n,\Ra^n)$. The cocycle super-rigidity
theorem then says that ``up to compact noise", the derivative of
any measure preserving $\G$ action on $\Ta^n$ looks measurably
like a constant linear map.  In fact, the cocycle superrigidity
theorems apply more generally to continuous actions on any
principal bundle $P$ over $M$ with fiber $H$, an algebraic group,
and in this context produces a measurable section
$s:M{\rightarrow}P$ satisfying equation $(\ref{equation:csr})$. So
in fact, cocycle superrigidity implies that for any action
preserving a finite measure on any manifold the derivative cocycle
looks measurably like a constant cocycle, up to compact noise.
That cocycle superrigidity provides information about actions of
groups on manifolds through the derivative cocycle  was first
observed in \cite{Fu}.  Zimmer originally proved cocycle
superrigidity in order to study orbit equivalence of group
actions.  For a recent survey of subsequent developments
concerning orbit equivalence rigidity and other forms of
superrigidity for cocycles, see \cite{Sl2}.

\subsection{Stability in hyperbolic dynamics.} \label{subsection:hyperbolic}

A diffeomorphism $f$ of a manifold $X$ is said to be {\em Anosov}
if there exists a continuous $f$ invariant splitting of the
tangent bundle $TX=E_{f}^u{\oplus}E_{f}^s$ and constants $a>1$ and
$C,C'>0$ such that for every $x{\in}X$,
\begin{enumerate}
\item $\|Df^n(v^u)\|{\geq}Ca^n\|v^u\|$ for all
$v^u{\in}E_{f}^u(x)$ and,

\item $\|Df^n(v^s)\|{\leq}C'a^{-n}\|v^s\|$ for all
$v^s{\in}E_{f}^s(x)$.
\end{enumerate}

\noindent We note that the constants $C$ and $C'$ depend on the
choice of metric, and that a metric can always be chosen so that
$C=C'=1$.  There is an analogous notion for a flow $f_t$, where
$TX=T\mathcal{O}{\oplus}E_{f_t}^u{\oplus}E_{f_t}^s$ where
$T{\mathcal{O}}$ is the tangent space to the flow direction and
vectors in $E_{f_t}^u$ (resp. $E_{f_t}^s$) are uniformly expanded
(resp. uniformly contracted) by the flow.  This notion was
introduced by Anosov and named after Anosov by Smale, who
popularized the notion in the United States \cite{An,Sm}. One of
the earliest results in the subject is Anosov's proof that Anosov
diffeomorphisms are {\em structurally stable}, or, in our language
$C^{1,1,1,0,0}$ locally rigid.  There is an analgous result for
flows, though this requires that one introduce a notion of time
change that we will not consider here.   Since Anosov also showed
that $C^2$ Anosov flows and diffeomorphisms are ergodic,
structural stability implies that the existence of an open set of
``chaotic" dynamical systems.

The notion of an Anosov diffeomorphism has had many interesting
generalizations, for example: Axiom A diffeomorphisms,
non-uniformly hyperbolic diffeomorphisms, and diffeomorphisms
admitting a dominated splitting.  The notion that has been most
useful in the study of local rigidity is the notion of a partially
hyperbolic diffeomorphism as introduced by Hirsch, Pugh and Shub.
Under strong enough hypotheses, these diffeomorphisms have a
weaker stability property similar to structural stability. More or
less, the diffeomorphisms are hyperbolic relative to some
foliation, and any nearby action is hyperbolic to some nearby
foliation. To describe more precisely the class of diffeomorphisms
we consider and the stability property they enjoy, we require some
definitions.

The use of the word {\em foliation} varies with context. Here a
{\em foliation by $C^k$ leaves} will be a continuous foliation
whose leaves are $C^k$ injectively immersed submanifolds that vary
continuously in the $C^k$ topology in the transverse direction. To
specify transverse regularity we will say that a foliation is
transversely $C^r$.  A foliation by $C^k$ leaves which is
tranversely $C^k$ is called simply a $C^k$ foliation. (Note our
language does not agree with that in the reference \cite{HPS}.)

Given an automorphism $f$ of a vector bundle $E{\rightarrow}X$ and
constants $a>b{\geq}1$, we say $f$ is {\em $(a,b)$-partially
hyperbolic} or simply {\em partially hyperbolic} if there is a
metric on $E$, a constant and $C{\geq}1$ and a continuous $f$
invariant, non-trivial splitting
$E=E_{f}^u{\oplus}E_{f}^c{\oplus}E_{f}^s$
 such that for every $x$ in $X$:
\begin{enumerate}

\item $\|f^n(v^u)\|{\geq}Ca^n\|v^u\|$ for all
$v^u{\in}E_{f}^u(x)$,

\item $\|f^n(v^s)\|{\leq}C{\inv}a^{-n}\|v^s\|$ for all
$v^s{\in}E_{f}^s(x)$ and

\item  $C{\inv}b^{-n}\|v^0\|<\|f^n(v^0)\|{\leq}C{b^n}\|v^0\|$ for
all $v^0{\in}E_{f}^c(x)$ and all integers $n$.

\end{enumerate}

\noindent A $C^1$ diffeomorphism $f$ of a manifold $X$ is {\em
$(a,b)$-partially hyperbolic} if the derivative action $Df$ is
$(a,b)$-partially hyperbolic on $TX$. We remark that for any
partially hyperbolic diffeomorphism, there always exists an {\it
adapted metric} for which $C=1$. Note that $E_{f}^c$ is called the
{\em central distribution} of $f$, $E_{f}^u$ is called the {\em
unstable distribution} of $f$ and $E_{f}^s$ the {\em stable
distribution} of $f$.

Integrability of various distributions for partially hyperbolic
dynamical systems is the subject of much research.  The stable and
unstable distributions are always tangent to invariant foliations
which we call the stable and unstable foliations and denote by
$\sw_f^s$ and $\sw_f^u$.  If the central distribution is tangent
to an $f$ invariant foliation, we call that foliation a {\em
central foliation} and denote it by $\sw^c_f$. If there is a
unique foliation tangent to the central distribution we call the
central distribution {\em uniquely integrable}. For smooth
distributions unique integrability is a consequence of
integrability, but the central distribution is usually not smooth.
If the central distribution of an $(a,b)$-partially hyperbolic
diffeomorphism $f$ is tangent to an invariant foliation $\sw^c_f$,
then we say $f$ is {\em $r$-normally hyperbolic to $\sw^c_f$} for
any $r$ such that $a>b^r$.  This is a special case of the
definition of $r$-normally hyperbolic given in \cite{HPS}.


Before stating a version of one of the main results of \cite{HPS},
we need one more definition.  Given a group $G$, a manifold $X$,
two foliations $\ff$ and $\ff'$ of $X$, and two actions $\rho$ and
$\rho'$ of $G$ on $X$, such that $\rho$ preserves $\ff$ and
$\rho'$ preserves $\ff'$, following \cite{HPS} we call $\rho$ and
$\rho'$ {\em leaf conjugate} if there is a homeomorphism $h$ of
$X$ such that:
\begin{enumerate}
\item $h(\ff)=\ff'$ and \item for every leaf $\fL$ of $\ff$ and
every $g{\in}G$, we have $h(\rho(g)\fL)=\rho'(g)h(\fL)$.
\end{enumerate}

\indent\indent
The map $h$ is then referred to as a {\em leaf conjugacy} between
$(X,\ff,\rho)$ and $(X,\ff',\rho')$.  This essentially means that
the actions are conjugate modulo the central foliations.

We state a special case of some the results of Hirsch-Pugh-Shub on
perturbations of partially hyperbolic actions of $\mathbb Z$, see
\cite{HPS}.  There are also analogous definitions and results for
flows.  As these are less important in the study of local
rigidity, we do not discuss them here.

\begin{theorem}
\label{theorem:hps} Let $f$ be an  $(a,b)$-partially hyperbolic
$C^k$ diffeomorphism of a compact manifold $M$ which is
$k$-normally hyperbolic to a $C^k$ central foliation $\sw^c_f$.
Then for any $\delta>0$, if $f'$ is a $C^k$ diffeomorphism of $M$
which is sufficiently $C^1$ close to $f$ we have the following:
\begin{enumerate}
\item $f'$ is $(a',b')$-partially hyperbolic, where
$|a-a'|<\delta$ and $|b-b'|<\delta$,  and  the splitting
$TM=E_{f'}^u{\oplus}E_{f'}^c{\oplus}E_{f'}^s$ for $f'$ is $C^0$
close to the splitting for $f$;

\item there exist $f'$ invariant foliations by $C^k$ leaves
$\sw^c_{f'}$ tangent to $E^c_{f'}$, which is close in the natural
topology on foliations by $C^k$ leaves to $\sw^c_f$,

\item there exists a (non-unique) homeomorphism $h$ of $M$ with
$h(\sw^c_f)=\sw^c_{f'}$, and $h$ is $C^k$ along leaves of
$\sw^c_f$, furthermore $h$ can be chosen to be $C^0$ small and
$C^k$ small along leaves of $\sw^c_{f}$

\item the homeomorphism $h$ is a leaf conjugacy between the
actions $(M,\sw^c_f,f)$ and $(M,\sw^c_{f'},f')$.
\end{enumerate}
\end{theorem}


\noindent Conclusion $(1)$ is easy and probably older than
\cite{HPS}. One motivation for Theorem \ref{theorem:hps} is to
study stability of dynamical properties of partially hyperbolic
diffeomorphisms. See the survey, \cite{BPSW}, for more discussion
of that and related issues.

\section{History}
\label{section:history}

In this section, we describe the history of the subject roughly to
the year 2000.  More recent developments will be discussed below.
Here we do not treat the subject entirely chronologically, but
break the discussion into four subjects: first, the study of local
rigidity of volume preserving actions, second the study of local
rigidity of certain non-volume preserving actions called boundary
actions, third the existence of (many) deformations of (many)
actions of groups that are typically quite rigid, and lastly a
brief discussion of infinitesimal rigidity. This is somewhat
ahistorical as the first results on smooth conjugacy of
perturbations of group actions appear in \cite{Gh1}, which we
describe in subsection \ref{subsection:boundaries}. While those
results are not precisely local rigidity results, they are clearly
related and the techniques involved inform some later approaches
to local rigidity.


\subsection{Volume preserving actions.}
\label{subsection:volumepreserving}

In this subsection, we discuss local rigidity of volume preserving
actions.  The acting groups will usually be lattices in higher
rank semi-simple Lie groups or higher rank free abelian groups.
Many of the results discussed here were motivated by conjectures
of Zimmer in \cite{Z4,Z5}.  The first result we mention, due to
Zimmer, does not prove local rigidity, but did motivate much later
work on the subject.

\begin{theorem}
\label{theorem:zimisom} Let $\G$ be a group with property $(T)$ of
Kazhdan and let $\rho$ be a Riemannian isometric action of $\G$ on
a compact manifold $M$.  Further assume the action is ergodic.
Then any $C^k$ action $\rho'$ which is $C^k$ close to $\rho$,
volume preserving and ergodic, preserves a $C^{k-3}$ Riemannian
metric.
\end{theorem}

\begin{remarks}
\begin{enumerate} \item Zimmer first
proved this theorem in \cite{Z1}, but only for $\G$ a lattice in a
semi-simple group, all of whose simple factors have real rank at
least two, and only producing a $C^0$ invariant metric for
$\rho'$. In \cite{Z2}, he gave the proof of the regularity stated
here and in \cite{Z4} he extended the theorem to all Kazhdan
groups.

\item In this theorem if $\rho'$ is $C^{\infty}$, the invariant
metric for $\rho'$ can also be chosen $C^{\infty}$.
\end{enumerate}
\end{remarks}

The first major result that actually produced a conjugacy between
the perturbed and unperturbed actions was due to Hurder,
\cite{H1}. Again, this result is not quite a local rigidity
theorem, but only a {\em deformation rigidity} theorem.  Hurder's
work is the first place where hyperbolic dynamics is used in the
theory, and is the beginning of a long development in which
hyperbolic dynamics play a key role.

\begin{theorem}
\label{theorem:hurder} The standard action of any finite index
subgroup of $SL(n,\Za)$ on the $n$ dimensional torus is
deformation rigid when $n{\geq}3$.
\end{theorem}

\noindent Hurder actually proves a much more general result.  His
theorem proves deformation rigidity of any group $\G$ acting on
the $n$ dimensional torus by linear transformations provided:
\begin{enumerate}
\item the set of periodic points for the $\G$ action is dense, and
\item the first cohomology of any finite index subgroup of $\G$ in
any $n$ dimensional representation vanishes, and \item the action
contains ``enough" Anosov elements.
\end{enumerate}\noindent Here we intentionally leave the meaning
of $(3)$ vague, as the precise notion needed by Hurder is quite
involved. To produce a continuous path of continuous conjugacies,
Hurder only need conditions $(1)$ and $(2)$ and the existence of
an Anosov diffeomorphism in the stabilizer of every periodic point
for the action. The additional Anosov elements needed in $(3)$ are
used to improve regularity of the conjugacy, and better techniques
for this are now available. Hurder's work has some applications to
actions of irreducible lattices in products of rank $1$ Lie
groups, which we discuss below in subsection
\ref{subsection:irredlattices}.  These applications do not appear
to be accessible by later techniques.

A key element in Hurder's argument is to use results of Stowe on
persistence of fixed points under perturbations of actions
\cite{St1,St2}. To use Stowe's result one requires that the
cohomology in the derivative representation at the fixed point
vanishes. This is where $(2)$ above is used.  Hurder constructs
his conjugacy by using the theorem of Anosov, mentioned above in
subsection \ref{subsection:hyperbolic}, that any Anosov
diffeomorphism is structurally stable.  This produces a conjugacy
$h$ for an Anosov element  $\rho(\gamma_0)$ which one then needs
to see is a conjugacy for the entire group action.  Hurder uses
Stowe's results to show that $h$ is conjugacy for the $\G$ actions
at all of the periodic points for the $\G$ action, and since
periodic points are dense it is then a conjugacy for the full
actions.  The precise argument using Stowe's theorem is quite
delicate and we do not attempt to summarize it here.  This
argument applies much more generally, see \cite[Theorem 2.9]{H1}.
That the conjugacy depends continuously, and in fact even
smoothly, on the original action is deduced from results on
hyperbolic dynamics in \cite{dlLMM}.

The first major development after Hurder's theorem was a theorem
of Katok and Lewis \cite{KL1} of which we state a special case:

\begin{theorem}
\label{theorem:katoklewis} Let $\G<SL(n,\Za)$ be a finite index
subgroup, $n>3$.  Then the linear action of $\G$ on $\Ta^n$ is
locally rigid.
\end{theorem}

It is worth noting that this theorem does not cover the case of
$n=3$.  A major ingredient in the proof is studying conjugacies
produced by hyperbolic dynamics for certain $\Za$ actions
generated by hyperbolic and partially hyperbolic diffeomorphisms
in both the original action and the perturbation.  The strategy is
to find a hyperbolic generating set and to show that the
conjugacies produced by the stability of those diffeomorphisms
agree.  A key ingredient idea is to show that they agree on the
set of periodic orbits.  Periodic orbits are then studied via
Theorem \ref{theorem:hps} for elements of $\G$ with large
centralizers and large central foliations. Periodic orbits are
detected as intersections of central foliations for different
elements, and this allows the authors to show that the structure
of the periodic set persists under deformations.  Elements with
large centralizers had previously been exploited by Lewis for
studying infinitesimal rigidity of similar actions. The authors
also exploit their methods to prove the following remarkable
result:

\begin{theorem}
\label{theorem:klabelian} Let $\Za^{n}$ be a maximal
diagonalizable (over $\Ra$) subgroup of $SL(n+1,\Za)$ where
$n{\geq}2$. The linear action of $\Za^{n}$ on $\Ta^{n+1}$ is
locally rigid.
\end{theorem}

This result is the first in a long series of results showing that
many actions of higher rank {\em abelian} groups are locally
rigid.  See Theorems \ref{theorem:ksabelian} and
\ref{theorem:damjanovichkatok} below for more instances of this
remarkable behavior.

The next major development occurs in a paper of Katok, Lewis and
Zimmer where Theorem \ref{theorem:katoklewis} is extended to cover
the case of $n=3$ as well as some more general groups acting on
tori.  Though this does not seem, on the face of it, to be a very
dramatic development, an important idea is introduced in this
paper.  The authors proceed by comparing the measurable data
coming from cocycle super-rigidity to the continuous data provided
by hyperbolic dynamics.  In this context, this essentially allows
the authors to show that the map $s$ in equation
$(\ref{equation:csr})$, described in the statement of cocycle
superrigidity given in subsection \ref{subsection:selbergweil}, is
continuous. This idea of comparing the output of cocycle
super-rigidity to information provided by hyperbolic dynamics has
played a major role in the development of both local and global
rigidity of group actions.

The results in the papers \cite{H1,KL1,KLZ} are all proven for
particular actions of particular groups, and in particular are all
proven for actions on tori.  The next sequence of developments was
a generalization of the ideas and methods contained in these
papers to fairly general Anosov actions of higher rank lattices on
nilmanifolds.  Part of this development takes place in the works
\cite{Q1,Q2,QY}.  A key difficulty in generalizing the early
approaches to rigidity of groups of toral automorphims is in
adapting the methods from hyperbolic dynamics which are used to
improve the regularity of the conjugacy.  In \cite{KS1,KS2}, Katok
and Spatzier developed a broadly applicable method for smoothing
conjugacies which depends on the theory of non-stationary normal
forms as developed by Guysinsky and Katok in \cite{GK,G}.  Two
main consequences of this method are:

\begin{theorem}
\label{theorem:ks} Let $G$ be a semi-simple Lie group with all
simple factors of real rank at least two, $\G<G$ a lattice, $N$ a
nilpotent Lie group and $\Lambda<N$ a lattice.  Then any affine
Anosov action of $\G$ on $N/{\Lambda}$ is locally rigid.
\end{theorem}

\noindent Here by {\em Anosov action}, we mean that some element
of $\G$ acts on $N/{\Lambda}$ as an Anosov diffeomorphism.

\begin{theorem}
\label{theorem:ksabelian} Let $N$ a nilpotent Lie group and
$\Lambda<N$ a lattice.  Let $\Za^d$ be a group of affine
transformations of $N/{\Lambda}$ such that the derivative action
(on a subgroup of finite index) is simultaneously diagonalizable
over $\Ra$ with no eigenvalues on the unit circle (i.e. on
subgroup of finite index, each element of $\Za^d$ is an Anosov
diffeomorphism which has semi-simple derivative). Then the $\Za^d$
action on $N/{\Lambda}$ is locally rigid.
\end{theorem}

\noindent Katok and Spatzier also apply their method to show that
for certain standard Anosov $\Ra^d$ actions the orbit foliation is
locally rigid.  I.e. any nearby action has conjugate orbit
foliation.  This yields interesting applications to rigidity of
boundary actions, see Theorem \ref{theorem:ksboundary} below.
Also, combined with other results of the same authors on rigidity
of cocycles over actions of Abelian groups, this yields local
rigidity of certain algebraic actions of $\Ra^d$, \cite{KS3,KS4}.
We state a special case of these results here:

\begin{theorem}
\label{theorem:katokspatzierRd} Let $G$ be an $\Ra$-split
semi-simple Lie group of real rank at least two.  Let $\Lambda<G$
be a cocompact lattice and let $\Ra^d<G$ be a maximal $\Ra$-split
subgroup.  Then the $\Ra^d$ action on $G/{\Lambda}$ is locally
rigid.
\end{theorem}

\begin{remarks}
\begin{enumerate} \item Here ``local
rigidity" has a slightly different meaning than above. Since the
automorphism group of $\Ra^d$ has non-trivial connected component,
it is possible to perturb the action by taking a small
automorphism of $\Ra^d$.  What is proven in this theorem is that
any small enough perturbation is conjugate to one obtained in this
way. \item Another approach to related cocycle rigidity results is
developed in the paper \cite{KNT}. \item The actual theorem in
\cite{KS2} is much more general.
\end{enumerate}
\end{remarks}

A key ingredient in the Katok--Spatzier method is to find
foliations which are orbits of transitive, isometric, smooth group
actions for both the perturbed and unperturbed action. To show
smoothness of the conjugacy, one constructs such group actions
that
\begin{enumerate} \item are intertwined by a continuous
conjugacy and \item exist on enough foliations to span all
directions in the space.
\end{enumerate} This proves that the conjugacy is ``smooth along
many directions" and one then uses a variety of analytic methods
to prove that the conjugacy is actually globally smooth.  The fact
that the transitive group exists and acts smoothly on the leaves
of some foliation for the unperturbed action is typically obvious.
One then uses the continuous conjugacy to define the group action
along leaves for the perturbed action and the fact that the
resulting action is smooth along leaves is verified using the
normal form theory. The foliations along which one builds
transitive group actions are typically central foliations for
certain special elements of the suspension of the action.  If the
original action is a $\Za^k$ action by automorphisms on some
nilmanifold $N/{\Lambda}$, the suspension of the action is the
left action of $\Ra^k$ on the solv-manifold
$M=(\Ra^k{\ltimes}N)/(\Za^k{\ltimes}\Lambda)$.  A typical one
parameter subgroup of $\Ra^k$ acts hyperbolically $M$, but certain
special directions in $\Ra^k$, those in so-called {\em Weyl
chamber walls} give rise to one parameter subgroups with
non-trivial central direction.  A key fact used in the argument is
that one can find another subgroup of $\Ra^k$ for which the
central foliation for some one parameter subgroup $\rho(t)$ is
also a dynamically defined, contracting foliation for the some
other element $a{\in}{\Ra^k}$.

All results quoted so far have strong assumptions on hyperbolicity
of the action.  For actions of semi-simple groups and their
lattices, the ultimate result on local rigidity in hyperbolic
context was proven by Margulis and Qian in \cite{MQ}.  This result
is for so-called weakly hyperbolic actions, which we define below.
This work proceeds by first using a comparison between hyperbolic
data and data from cocycle superrigidity to produce a continuous
conjugacy between $C^1$ close actions and then uses an adaptation
of the Katok--Spatzier smoothing method mentioned above.  A key
technical innovation in this work is the choice of cocycle to
which cocycle super-rigidity is applied.  In all work to this
point, it was applied to the derivative cocycle.  Here it is
applied to a cocycle that measures the difference between the
action and the perturbation.  To illustrate the idea, we give the
definition of this cocycle, which we refer to as the {\em
Margulis--Qian cocycle}, in the special case of actions by left
translations. As this construction is quite general, we will let
$D$ be the acting group.  Let the $D$ action $\rho$ on
$H/{\Lambda}$ be defined via a homomorphism
$\pi_0:D{\rightarrow}H$. Let $\rho'$ be a perturbation of $\rho$.
If $D$ is connected it is clear that the action lifts to $\tilde
H$ and therefore to $H$. If $D$ is discrete, this lifting still
occurs, since the obstacle to lifting is a cohomology class in
$H^2(D, \pi_1(H/{\Lambda}))$ which  does not change under a small
perturbation of the action. (A direct justification without
reference to group cohomology can be found in \cite{MQ} section
2.3.)  Write the lifted actions of $D$ on $H$ by $\tilde \rho$ and
${\tilde \rho}'$ respectively. We can now define a cocycle
$\alpha:D{\times}H{\rightarrow}H$ by
$${\tilde \rho}'(g)x={\alpha}(g,x)x$$
for any $g$ in $D$ and any $x$ in $H$. It is easy to check that
this is a cocycle and that it is right $\Lambda$ invariant, and so
defines a cocycle $\alpha:D{\times}H/{\Lambda}{\rightarrow}H$. See
\cite{MQ} section 2 or \cite{FM1} section 6 for more discussion as
well as for more general variants on this definition.  We remark
that the use of this cocycle allowed Margulis and Qian to prove
the first local rigidity results for volume preserving actions of
lattices that have no global fixed point.  The construction of
this cocycle is inspired by a cocycle used by Margulis in his
first proof of superrigidity.  This is the cocycle
$G{\times}G/\G{\rightarrow}\G$ defined by the choice of a
fundamental domain for $\G$ in $G$.  See Example $4$ in subsection
\ref{subsection:irredlattices} for a more explicit description.

The work of Margulis and Qian applies to actions which satisfy the
following condition.  This condition essentially says that the
action is hyperbolic in all possible directions, at least for some
element of the acting group.  It is easy to construct weakly
hyperbolic actions of lattices, in particular $\G$ acting on
$G/\Lambda$ where $G$ is a simple Lie group  and $\Lambda<G$ is a
cocompact lattice and $\Gamma$ is any other lattice in $G$.  It is
important to note that for this example, no element acts as an
Anosov diffeomorphism and, with an appropriate choice of $\G$ and
$\Lambda$, there are no finite $\G$ orbits.

\begin{defn}
\label{defn:weaklyhyperbolic} An action $\rho$ of a group $D$ on a
manifold $M$ is called {\em weakly hyperbolic} if there exist
elements $d_1,{\ldots},d_k$ and constants $a_i>b_i{\geq}1$ for
$i=1,\ldots,k$ such that each $\rho(d_i)$ is $(a_i,b_i)$-partially
hyperbolic in the sense of subsection \ref{subsection:hyperbolic}
and we have $TM=\sum E^s_{\rho(d_i)}$. I.e. there are partially
hyperbolic elements whose stable (or unstable) directions span the
tangent space at any point.
\end{defn}

\begin{theorem}
\label{theorem:margulisqian} Let $H$ be a real algebraic Lie group
and $\Lambda<H$ a cocompact lattice.  Assume $G$ is a semisimple
Lie group with all simple factors of real rank at least two and
$\G<G$ is a lattice.  Then any weakly hyperbolic, affine algebraic
action of $\G$ or $G$ on $H/{\Lambda}$ is locally rigid.
\end{theorem}

\begin{remark}
This is somewhat more general than the
result claimed in \cite{MQ}, as they only work with certain
special classes of affine actions, which they call standard. This
result can proven by the methods of \cite{MQ}, and a proof in
precisely this generality can be read out of \cite{FM1,FM3},
simply by assuming that the common central direction for the
acting group is trivial.
\end{remark}

The next major result was a remarkable theorem of Benveniste
concerning isometric actions.  This is a stronger result than
Theorem \ref{theorem:zimisom} because it actually produces a
conjugacy, but is weaker in that it requires much stronger
assumptions on the acting group.

\begin{theorem}
\label{theorem:benveniste} Let $\G$ be a cocompact lattice in a
semi-simple Lie group with all simple factors of real rank at
least two.   Let $\rho$ be an isometric $\G$ action on a compact
manifold $M$.  Then $\rho$ is locally rigid.
\end{theorem}

\noindent The proof of this theorem is inspired by the work of
Calabi, Vesentini and Weil in the original proofs of Theorem
\ref{theorem:rigidityoflattices} and is based on showing that
certain deformations of foliated geometric structures are trivial.
The argument is much more difficult than the classical case and
uses Hamilton's implicit function theorem.  This is the first
occasion on which analytic methods like KAM theory or hard
implicit function theorems appear in work on local rigidity of
group actions.  More recently these kinds of methods have been
applied more systematically, see subsection \ref{subsection:kam}.

The theorems described so far concern actions that are either
isometric or weakly hyperbolic. There are many affine actions
which satisfy neither of these dynamical hypotheses, but are
genuinely partially hyperbolic.  Local rigidity results for
actions of this kind first arise in work of Nitica and Torok. We
state special cases of two of their theorems:

\begin{theorem}
\label{theorem:niticatorok1} Let $\G<SL(n,\Za)$ be a finite index
subgroup with $n{\geq}3$.  Let $\rho_1$ be the standard $\G$
action on $\Ta^n$ and let $\rho$ be the diagonal $\G$ action on
$\Ta^n{\times}\Ta^m$ defined by $\rho_1$ on the first factor and
the trivial action on the second factor.  The action $\rho$ is
deformation rigid.
\end{theorem}

\begin{theorem}
\label{theorem:niticatorok2} Let $\G,\rho_1,\rho$ be as above and
further assume that $m=1$.  The action $\rho$ is locally rigid.
\end{theorem}

\begin{remarks}
\begin{enumerate}
\item Nitica and Torok prove more general theorems in which both
$\G$ and $\rho_1$ can be more general.  The exact hypotheses
required are different in the two theorems.

\item We are being somewhat ahistorical here, Theorem
\ref{theorem:niticatorok1} predates the work of Margulis and Qian.

\item The conjugacy produced in the papers \cite{NT1,NT2,T} is
never $C^{\infty}$, but only $C^k$ for some choice of $k$.  The
choice of $k$ is essentially free and determines the size of
perturbations or deformations that can be considered.  It should
be possible to produce a $C^{\infty}$ conjugacy by combining the
arguments in these papers with arguments in \cite{FM2,FM3}, see
the end of subsection \ref{subsection:fmsketch} for some
discussion.
\end{enumerate}
\end{remarks}

The work of Nitica and Torok is quite complex, using several
different ideas.  The most novel is to study rigidity of cocycles
over hyperbolic dynamical systems taking values in diffeomorphism
groups.  The dynamical system is either the action $\rho_1$ in
Theorem \ref{theorem:niticatorok1} or \ref{theorem:niticatorok2}
or it's restriction to {\em any} sufficiently generic subgroup
containing an Anosov diffeomorphism of $\Ta^n$, and the target
group is the group of diffeomorphisms of $\Ta^m$. This part of the
work is inspired by a classical theorem of Livsic and the proof
his modelled on his proof.  To reduce the rigidity question to the
cocycle question is quite difficult and depends on an adaptation
of the work of \cite{HPS} discussed in subsection
\ref{subsection:hyperbolic} as well as use of results of Stowe
\cite{St1,St2}. Regrettably, the technology seems to limit the
applicability of the ideas to diagonal actions $\rho_1{\times}\Id$
on products $M{\times}N$ where the action on $M$ has many periodic
points and the action on $N$ is trivial. Theorem
\ref{theorem:niticatorok2} also depends on the acting group having
property $(T)$ of Kazhdan.  The method of proof of Theorem
\ref{theorem:niticatorok1} has additional applications, see
particularly subsection \ref{subsection:irredlattices} below.

To close this section, we remark that local rigidity is often
considerably easier in the analytic setting.  Not much work has
been done in this direction, but there is an interesting note of
Zeghib \cite{Zg}. A sample result is the following:

\begin{theorem}
\label{theorem:zeghib} Let $\G<SL(n,\Za)$ be a subgroup of finite
index and let $\rho$ be the standard action of $\G$ on $\Ta^n$.
Then any analytic action close enough to $\rho$ is analytically
conjugate to $\rho$.  Furthermore, if $M$ is a compact analytic
manifold on which $\G$ acts trivially and we let $\tilde \rho$ be
the diagonal action of $\G$ on $\Ta^n{\times}M$, then $\tilde
\rho$ is also locally rigid in the analytic category.
\end{theorem}

\noindent Zeghib also proves a number of other interesting results
for both volume preserving and non-volume preserving actions and
it is clear that his method has applications not stated in his
note.  The key point for all of his arguments is a theorem of Ghys
and Cairns that says that one can linearize an analytic action of
a higher rank lattice in a neighborhood of a fixed point.  Zeghib
proves his results by using results of Stowe \cite{St1,St2} to
find fixed points for the perturbed action and then studying the
largest possible set to which the linearization around this point
can be extended. We end this section by stating the theorem of
Cairns and Ghys from \cite{CGh} which Zeghib uses.

\begin{theorem}
\label{theorem:ghyscairns} Let $G$ be a semisimple Lie group of
real rank at least two with no compact factors and finite center
and let $\G<G$ be a lattice.  Then every analytic action of $\G$
with a fixed point $p$ is analytically linearizable in a
neighborhood of $p$.
\end{theorem}

\begin{remarks} \begin{enumerate} \item By analytically
linearizable in a neighborhood of $p$, we mean that there exists a
neighborhood $U$ of $p$ and an analytic diffeomorphism $\phi$ of
$U$ into the ambient manifold $M$ such that the action of $\G$,
conjugated by $\phi$ is the restriction of a linear action to
$\phi(U)$. \item In the same paper, Ghys and Cairns give an
example of a $C^{\infty}$ action of $SL(3,\Za)$ on $\Ra^8$ fixing
the origin, which is not $C^0$ linearizable in any neighborhood of
the origin. So the assumption of analyticity in the theorem is
necessary.
\end{enumerate}
\end{remarks}

\subsection{Actions on boundaries.}
\label{subsection:boundaries}

In this subsection we discuss rigidity results for groups acting
on homogeneous spaces known as ``boundaries".  In contrast to the
last section, the actions we describe here never preserve a volume
form, or even a Borel measure. We will not discuss here all the
geometric, function theoretic or probabalistic reasons why these
spaces are termed boundaries, but merely describe examples. For
us, if $G$ is a semisimple Lie group, then a {\em boundary} of $G$
is a space of the form $G/P$ where $P$ is a connected Lie subgroup
of $G$ such that the quotient $G/P$ is compact.  The groups $P$
having this property are often called {\em parabolic} subgroups.
The space $G/P$ is also considered a boundary for any lattice
$\G<G$. For more precise motivation for this terminology, see
\cite{Fu1,Fu2,Mo1,Ma3}.

The simplest example of a boundary is for the group $SL(2,\Ra)$ in
which case the only choice of $P$ resulting in a non-trivial
boundary is the group of upper (or lower) triangular matrices. The
resulting quotient is naturally diffeomorphic to the circle and
$SL(2,\Ra)$ acts on this circle by the action on rays through the
origin in $\Ra^2$.  We can restrict this $SL(2,\Ra)$ action to any
lattice $\G$ in $SL(2,\Ra)$.  The following remarkable theorem was
first proved by Ghys in \cite{Gh1}:

\begin{theorem}
\label{theorem:ghys} Let $\G<SL(2,\Ra)$ be a cocompact lattice and
let $\rho$ be the action of $\G$ on $S^1$ described above. If
$\rho'$ is any perturbation of $\rho$, then $\rho'$ is smoothly
conjugate to an action defined by another embedding $\pi'$ of $\G$
in $SL(2,\Ra)$ close to the original embedding.  In particular
$\pi'(\G)$ is a cocompact lattice in $SL(2,\Ra)$.
\end{theorem}

\noindent Ghys gives two proofs of this fact, one in \cite{Gh1}
and another different one in \cite{Gh2}.  A third and also
different proof is in later work of Kononenko and Yue \cite{KY}.
Ghys' first proof derives from a remarkable global rigidity result
for actions on certain three dimensional manifolds by the affine
group of the line, while his second derives from rigidity results
concerning certain Anosov flows on three dimensional manifolds. We
remark that the fact that $\rho'$ is continuously conjugate to an
action defined by a nearby embedding into $SL(2,\Ra)$  was known
and so Theorem \ref{theorem:ghys} can be viewed as a regularity
theorem though this is not how the proof proceeds.

Both of Ghys' proofs pass through a statement concerning local
rigidity of foliations.  This uses the following variant on the
construction of the induced action.  Let $\G$ be a cocompact
lattice in $SL(2,\Ra)$, and let $\rho$ be the $\G$ action on $S^1$
defined by the action of $SL(2,\Ra)$ there.  There is also a $\G$
action on the hyperbolic plane $\Ha^2$.  We form the manifold
$(\Ha^2{\times}S^1)/\G_{\rho}$, where $\G$ acts diagonally. This
manifold is diffeomorphic to the unit tangent bundle of $\Ha^2/\G$
which is also diffeomorphic to $SL(2,\Ra)/{\G}$, and the foliation
by planes of the form $\Ha^2{\times}\{point\}$ is the weak stable
foliation for the geodesic flow and also the orbit foliation for
the action of the affine group. Given a $C^r$ perturbation $\rho'$
of the $\G$ action on $S^1$, we can form the corresponding bundle
$(\Ha^2{\times}S^1)_{\rho'}/\G_{\rho'}$, and the foliation by
planes of the form $\Ha^2{\times}\{point\}$ is $C^r$ close to
analogous foliation in $(\Ha^2{\times}S^1)\G_{\rho}$. To show that
$\rho$ and $\rho'$ are conjugate, it suffices to find a
diffeomorphism of $(\Ha^2{\times}S^1)/\G_{\rho}$ conjugating the
two foliations. Both of Ghys' proofs proceed by constructing such
a conjugacy of foliations.  This reduction to studying local
rigidity of foliations has further applications in slightly
different settings, see Theorems \ref{theorem:ksboundary} and
\ref{theorem:kanai1} below.

In later work, Ghys proved a remarkable result which characterized
an entire connected component of the space of actions of $\G$ on
$S^1$.  Let $X$ be the component of $\Hom(\G,\Diff(S^1))$
containing the actions described in Theorem \ref{theorem:ghys}. In
\cite{Gh3}, Ghys showed that this component consisted entirely of
actions conjugate to actions defined by embeddings $\pi'$ of $\G$
into $SL(2,\Ra)$ where $\pi'(\G)$ is a cocompact lattice.  This
result builds on earlier work of Ghys where a similar result was
proven concerning $\Hom(\G,\Homeo(S^1))$.  A key ingredient is the
use of the Euler class of the action, viewed as a bounded cocycle.

As mentioned above, Ghys' method of reducing local rigidity of an
action to local rigidity of a foliation has had two more
applications.  The first of these is due to Katok and Spatzier
\cite{KS2}.

\begin{theorem}
\label{theorem:ksboundary} Let $G$ be a semisimple Lie group with
no compact factors and real rank at least two.  Let $\G<G$ be a
cocompact lattice and $B=G/P$ a boundary for $\G$.  Then the $\G$
action on $G/P$ is locally rigid.
\end{theorem}

\noindent The proof of this result uses an argument similar to
Ghys' to reduce to a need to study regularity of foliations for
perturbations of the action of certain connected abelian subgroups
of $G$ on $G/{\Gamma}$. The result used in the proof here is the
same as the one used in the proof of Theorem
\ref{theorem:katokspatzierRd}. For $\G$ a lattice but not
cocompact, some partial results are obtained by Yaskolko in his
Ph.D. thesis \cite{Yk}.

Following a similar outline, Kanai proved the following:

\begin{theorem}
\label{theorem:kanai1} Let $G=SO(n,1)$ and $\G<G$ be a cocompact
lattice.  Then the action of $\G$ on the boundary $G/P$ is locally
rigid.
\end{theorem}

\noindent Partial results in this direction were proven earlier by
Chengbo Yue.  Yue also proves partial results in the case where
$SO(1,n)$ is replaced by any rank $1$ non-compact simple Lie
group.

In somewhat earlier work, Kanai had also proven a special case of
Theorem \ref{theorem:ksboundary}.  More precisely:

\begin{theorem}
\label{theorem:kanai2} Let $\G<SL(n,\Ra)$ be a cocompact lattice
where $n{\geq}21$ and let $\rho$ be the $\G$ action on $S^{n-1}$
by acting on the space of rays in $\Ra^n$.  Then $\rho$ is locally
rigid.
\end{theorem}

\noindent  Kanai's proof proceeds in two steps.  In the first
step, he uses Thomas' notion of a projective connection to reduce
the question to one concerning vanishing of certain cohomology
groups.  In the second step, he uses stochastic calculus to prove
a vanishing theorem for the relevant cohomology groups. The first
step is rather special, and is what restricts Kanai's attention to
spheres, rather than other boundaries, which are Grassmanians. The
method in the second step seems a good deal more general and
should have further applications, perhaps in the context of
Theorem \ref{theorem:affine} below.  While the approach here is
similar in spirit to the work of Benveniste in \cite{Be1}, it
should be noted that Kanai does not use a hard implicit function
theorem.

We end this subsection by recalling a construction due to Stuck,
which shows that much less rigidity should be expected from
actions which do not preserve volume \cite{Sk}.  Let $G$ be a
semi-simple Lie group and $P<G$ a minimal parabolic.  Then there
always exists a homomorphism $\pi:P{\rightarrow}\Ra$.  Given any
manifold $M$ and any action $s$ of $\Ra$ on $M$, we can then form
the induced $G$ action $\rho_s$ on $(G{\times}M)/P$ where $P$ acts
on $G$ on the left and on $M$ by $\pi{\circ}s$.  Varying the
action $s$ varies the action $\rho_s$.  It is easy to see that if
$\rho_s$ and $\rho_{s'}$ are two such actions, then they are
conjugate as $G$ actions if and only if $s$ and $s'$ are
conjugate.  If one picks an irreducible lattice $\G$ in
$G{\times}G$, project $\G$ to $G$ and restricts the actions to
$\G$, then it is also easy to see that the restriction of $\rho_s$
and $\rho_{s'}$ to $\G$ are conjugate if and only if $s$ and $s'$
are conjugate.  The author does not know a proof that this is also
true if one simply takes a lattice $\G$ in $G$, but believes that
this is also true and may even be known.


\subsection{``Flexible" actions of rigid groups.}
\label{subsection:flexible}

In this subsection, I discuss a sequence of results concerning
flexible actions of large groups.  More or less, the sequence of
examples provides counter-examples to most naive conjectures of
the form ``all of actions of some lattice $\G$ are locally rigid."
There are some groups, for example compact groups and finite
groups, all of whose smooth actions are locally rigid.  It seems
likely that there should be infinite discrete groups with this
property as well, but the constructions in this subsection show
that one must look beyond lattices in Lie groups for examples.

Essentially all of the examples given here derive from the simple
construction of ``blowing up" a point or a closed orbit, which was
introduced to this subject in \cite{KL2}.  The further
developments after that result are all further elaborations on one
basic construction.  The idea is to use the ``blow up"
construction to introduce distinguished closed invariant sets
which can be varied in some manner to produce deformations of the
action.  The ``blow up" construction is a classical tool from
algebraic geometry which takes a manifold $N$ and a point $p$ and
constructs from it a new manifold $N'$ by replacing $p$ by the
space of directions at $p$. Let ${\Ra}P^l$ be the $l$ dimensional
projective space.  To blow up a point, we take the product of
$N{\times}{\Ra}P^{\dim(N)}$ and then find a submanifold where the
projection to $N$ is a diffeomorphism off of $p$ and the fiber of
the projection over $p$ is ${\Ra}P^{\dim(N)}$.   For detailed
discussion of this construction we refer the reader to any
reasonable book on algebraic geometry.

The easiest example to consider is to take the action of
$SL(n,\Za)$, or any subgroup  $\G<SL(n,\Za)$ on the torus $\Ta^n$
and blow up the fixed point, in this case the equivalence class of
the origin in $\Ra^n$. Call the resulting manifold $M$. Provided
$\G$ is large enough, e.g. Zariski dense in $SL(n,\Ra)$,  this
action of $\G$ does not preserve the measure defined by any volume
form on $M$. A clever construction introduced in \cite{KL2} shows
that one can alter the standard blowing up procedure in order to
produce a one parameter family of $SL(n,\Za)$ actions on $M$, only
one of which preserves a volume form. This immediately shows that
this action on $M$ admits perturbations, since it cannot be
conjugate to the nearby, non-volume preserving actions.
Essentially, one constructs different differentiable structures on
$M$ which are diffeomorphic but not equivariantly diffeomorphic.

After noticing this construction, one can proceed to build more
complicated examples by passing to a subgroup of finite index, and
then blowing up several fixed points.  One can also glue together
the ``blown up" fixed points to obtain an action on a manifold
with more complicated topology. See \cite{KL2,FW} for discussion
of the topological complications one can introduce.

In \cite{Be2} it is observed that a similar construction can be
used for the action of a simple group $G$ by left translations on
a homogeneous space $H/\Lambda$ where $H$ is a Lie group
containing $G$ and $\Lambda<H$ is a cocompact lattice.  Here we
use a slightly more involved construction from algebraic geometry,
and ``blow up" the directions normal to a closed submanifold. I.e.
we replace some closed submanifold $N$ in $H/{\Lambda}$ by the
projectived normal bundle to $N$.  In all cases we consider here,
this normal bundle is trivial and so is just $N{\times}{\Ra}P^{l}$
where $l=\dim(H)-\dim(N)$.

Benveniste used his construction to produce more interesting
perturbations of actions of higher rank simple Lie group $G$ or a
lattice $\G$ in $G$.  In particular, he produced volume preserving
actions which admit volume preserving perturbations.  He does this
by choosing $G<H$ such that not only are there closed $G$ orbits
but so that the centralizer $Z=Z_H(G)$ of $G$ in $H$  has
no-trivial connected component.  If we take a closed $G$ orbit
$N$, then any translate $zN$ for $z$ in $Z$ is also closed and so
we have a continuum of closed $G$ orbits.  Benveniste shows that
if we choose two closed orbits $N$ and $zN$ to blow up and glue,
and then vary $z$ in a small open set, the resulting actions can
only be conjugate for a countable set of choices of $z$.

This construction is further elaborated in \cite{F1}. Benveniste's
construction is not optimal in several senses, nor is his proof of
rigidity.  In \cite{F1}, I give a modification of the construction
that produces non-conjugate actions for every choice of $z$ in a
small enough neighborhood.  By blowing up and gluing more pairs of
closed orbits, this allows me to produce actions where the space
of deformations contains a submanifold of arbitrarily high, finite
dimension. Further, Benveniste's proof that the deformation are
non-trivial is quite involved and applies only to higher rank
groups. In \cite{F1}, I give a different proof of non-triviality
of the deformations, using consequences of Ratner's theorem to due
Witte and Shah \cite{R,Sh,W}. This shows that the construction
produces non-trivial perturbations for any semisimple $G$ and any
lattice $\G$ in $G$.

In \cite{BF} we show that none of these actions preserve any rigid
geometric structure in the sense of Gromov.  It is possible that
any action of a higher rank lattice which preserves a rigid
geometric structure is locally rigid. It is also possible that any
such action is generalized quasi-affine.

\subsection{Infinitesimal rigidity.}
\label{subsection:infrigid}

In \cite{Z4}, Zimmer introduced a notion of infinitesimal rigidity
motivated by Weil's Theorem \ref{theorem:weil} and the analogy
between finite dimensional Lie algebras and vector fields. Let
$\rho$ be a smooth action of a group $\G$ on a manifold $M$, then
$\rho$ is {\em infinitesimally rigid} if
$H^1(\G,\Vect^{\infty}(M))=0$.  Here the $\G$ action on
$\Vect^{\infty}(M)$ is given by the derivative of $\rho$. The
notion of infinitesimal rigidity was introduced with the hope that
one could prove an analogue of Weil's Theorem \ref{theorem:weil}
and then results concerning infinitesimal rigidity would imply
results concerning local rigidity.  Many infinitesimal rigidity
results were then proven, see \cite{H2,Ko,L,LZ,Q3,Z6}. For some
more results on infinitesimal rigidity, see Theorems
\ref{theorem:irredlattices} and \ref{theorem:clozel}.  Also see
subsection \ref{subsection:kam} for a discussion of known results
on the relation between infinitesimal and local rigidity.

\section{Recent developments}
\label{section:recent}

In this section we discuss the most recent dramatic developments
in the field.  The first subsection discusses work of the author
and Margulis on rigidity of actions of higher rank groups and
lattices. Our main result is that if $H$ is the real points of an
algebraic group defined over $\Ra$ and $\Lambda<H$ is a cocompact
lattice, then any affine action of $G$ or $\G$ on $H/{\Lambda}$ is
locally rigid. This work is quite involved and spans a sequence of
three long papers \cite{FM1,FM2,FM3}.  One of the main goals of
subsection \ref{subsection:fmsketch} is to provide something of a
``reader's guide" to those papers.

The second subsection discusses some recent developments involving
more geometric and analytic approaches to local rigidity. Till
this point, the study of local rigidity of group actions has been
dominated by algebraic ideas and hyperbolic dynamics with the
exception of \cite{Ka1} and \cite{Be2}.  The results described in
subsection \ref{subsection:kam} represent (the beginning of) a
dramatic development in analytic and geometric techniques. The
first of these is the work of Damjanovich and Katok on local
rigidity of certain partially hyperbolic affine actions of abelian
groups on tori using a KAM approach \cite{DK1,DK2}.  The second is
the author's proof of a criterion for local rigidity of groups
actions modelled on Weil's Theorem \ref{theorem:weil} and proven
using Hamilton's implicit function theorem \cite{F2}.  This result
currently has an unfortunate ``side condition" on second
cohomology that makes it difficult to apply.

The final subsection concerns a few other very recent results and
developments that the author feels point towards the future
development of the field.

\subsection{The work of Margulis and the author.}
\label{subsection:fmsketch}

Let $G$ be a (connected) semi-simple Lie group with all simple
factors of real rank at least two, and $\Gamma<G$ is a lattice.
The main result of the papers \cite{FM1,FM2,FM3} is:

\begin{theorem}
\label{theorem:main} Let $\rho$ be a volume preserving
quasi-affine action of $G$ or $\Gamma$ on a compact manifold $X$.
Then the action locally rigid.
\end{theorem}

\begin{remarks}
\begin{enumerate}
\item This result subsumes essentially all of the theorems in
subsection \ref{subsection:volumepreserving}, excepting those
concerning actions of abelian groups.

\item In \cite{FM3} we also achieve some remarkable results for
perturbations of very low regularity. In particular, we prove that
any perturbation which is a $C^3$ close $C^3$ action is conjugate
back to the original action by a $C^2$ diffeomorphism. \item The
statement here is slightly different than that in \cite{FM3}. Here
$X=H/L{\times}M$ with $L$ cocompact, while there
$X=H/{\Lambda}{\times}M$ with $\Lambda$ discrete and cocompact. An
essentially algebraic argument using results in \cite{W1}, shows
that possibly after changing $H$ and $M$, these hypotheses are
equivalent.
\end{enumerate}
\end{remarks}

Another main result of the research resulting in Theorem
\ref{theorem:main} is the following:

\begin{theorem}
\label{theorem:isomrigid} Let $\Gamma$ be a discrete group with
property $(T)$.  Let $X$ be a compact smooth manifold, and let
$\rho$ be a smooth action of $\Gamma$ on $X$ by Riemannian
isometries. Then $\rho$ is locally rigid.
\end{theorem}

\begin{remarks}
\begin{enumerate} \item A key step in the
proof of Theorem \ref{theorem:main} is a foliated version of
Theorem \ref{theorem:isomrigid}. \item As in Theorem
\ref{theorem:main}, there is a finite regularity version of
Theorem \ref{theorem:isomrigid} and it's foliated generalization,
we refer the reader to \cite{FM2} for details.
\end{enumerate}
\end{remarks}

The remainder of this subsection will consist of a  sketch of the
proof of Theorem \ref{theorem:main}.  The intention is essentially
to provide a reader's guide to the three papers
\cite{FM1,FM2,FM3}.  Throughout the remainder of this subsection,
to simplify notation, we will discuss only the case of affine $\G$
actions on $H/\Lambda$ with $\Lambda$ a cocompact lattice. The
proof for connected groups and quasi-affine action on
$X=H/{\Lambda}{\times}M$ is similar. To further simplify the
discussion, we assume that $\rho'$ is a $C^{\infty}$ perturbation
of $\rho$.

\subsection*{\bf Step 1: An invariant ``central" foliation for the
perturbed action and leaf conjugacy.}

To begin the discussion and the proof, we need some knowledge of
the structure of the affine actions considered. By \cite[Theorem
6.4]{FM1}, there is a finite index subgroup $\G'<\G$ such that the
action of $\G'$ on $H/{\Lambda}$ is given by a homomorphism
$\sigma:\G'{\rightarrow}\Aut(H){\ltimes}H$. We simplify the
discussion by assuming $\G=\G'$ throughout.  Using Margulis
superrigidity theorems, which are also used in the proof of
\cite[Theorem 6.4]{FM1}, it is relatively easy to understand the
maximal central foliation $\ff$ for $\rho$:  there is a subgroup
$Z<H$ whose orbit foliation is exactly the central foliation. For
example, if $G<H$ acts on $H/{\Lambda}$ by left translations and
$\rho$ is restriction of that action to $\G$, then $Z=Z_H(G)$. For
details on what $Z$ is more generally, see \cite{FM1}.

Given a perturbation $\rho'$ of $\rho$, we begin by finding a
$\rho'$ invariant foliation $\ff'$ and a leaf conjugacy $\phi$
from  $(H/\Lambda,\rho,\ff)$ to $(H/\Lambda,\rho',\ff')$. To do
this, we apply a result concerning local rigidity of cocycles over
actions of higher rank groups and lattices to the Margulis--Qian
cocycle defined by the perturbation. As the statements of the
local rigidity results for cocycles are somewhat technical, we
refer the reader to \cite[Theorems 1.1 and 5.1]{FM1}.  Those
Theorems are proven in Section $5$ of that paper using results in
Section $4$ concerning orbits in representation varieties as well
as the cocycle superrigidity theorems.  The construction of the
leaf conjugacy is completed in \cite[Section 2.2]{FM3} using
\cite[Theorem 1.8]{FM1}.  We remark that we actually construct
$\phi{\inv}$ rather than $\phi$.  The paper \cite{FM1} also
contains a proof of superrigidity for cocycles that results in
many technical improvements to that result.

\subsection*{\bf Step 2: Smoothness of the central foliation,
reduction to a foliated perturbation.}

The next step in the proof is to show that $\ff'$ is a foliation
by smooth leaves.  In fact, it is only possible to show at this
point that it is a foliation by $C^k$ leaves for some $k$
depending on $\rho$ and $\rho'$ and particularly on the $C^1$ size
of the perturbation. This is done using the work of Hirsch, Pugh
and Shub described in subsection \ref{subsection:hyperbolic}.  If
the central foliation for $\rho$ is the central foliation for
$\rho(\g)$ for some single element $\g$ in $\G$, this amounts to
showing that the foliation $\ff'$ constructed in step one is the
same foliation as the central foliation for $\rho'(\g)$
constructed in the proof of Theorem \ref{theorem:hps}.  To prove
this, one needs to analyze the proof of Theorem \ref{theorem:hps}.
More generally, we show, in \cite[Section 3.2]{FM3} that there
is a finite collection of elements $\g_1,\ldots,\g_k$ in $\G$ such
that each leaf of the foliation $\ff$ is a transverse intersection
of central leaves of $\rho(\g_1),\ldots,\rho(\g_k)$. One then
needs to combine an analysis of the proof of Theorem
\ref{theorem:hps} with some arguments concerning persistence of
transversality under certain kinds of perturbations.  This
argument is carried out in \cite[Section 3.3]{FM3}.

Once we know that $\ff'$ is $C^k$, it is easy to see that the leaf
conjugacy $\phi$ is $C^k$ and $C^k$ small along leaves of $\ff$
though all derivatives are only continuous in the transverse
direction.  Conjugating the action $\rho'$ by $\phi$, we obtain a
new action $\rho''$ on $H/{\Lambda}$ which preserves $\ff$.  This
action is only continuous, but it is $C^0$ close to $\rho$ and
$C^k$ and $C^k$ close to $\rho$ along leaves of $\ff$.  In
\cite{FM2,FM3} we refer to perturbations of this type as {\em
foliated perturbations}.

\subsection*{\bf Step 3: Conjugacy along the central foliation.}

The next step is to apply a foliated generalization of Theorem
\ref{theorem:isomrigid} to the actions $\rho$ and $\rho''$. The
exact result we apply is \cite[Theorem 2.11]{FM2} which produces a
semi-conjugacy $\psi$ between $\rho$ and $\rho'$.  This result is
somewhat involved to state and the regularity of $\psi$ is hard to
describe.  The map $\psi$ is $C^{k-1-\varepsilon}$ along the
leaves of $\ff$ at almost all points in $H/{\Lambda}$, for
$\varepsilon$ depending on the size of the perturbation, but only
transversely measurable. In addition, the map $\psi$ satisfies a
certain Sobolev estimate, that implies that it is
$C^{k-1-\varepsilon}$ small in a small ball in $\ff$ at most
points, and that the $C^{k-1-\varepsilon}$ norm is only large on
very small sets. Rather than try to make this precise here, we
include a sketch of the proof of Theorem \ref{theorem:isomrigid}.
Before doing so we remark that the map $\varphi=\phi{\circ}\psi$
is a semiconjugacy from between the $\G$ action $\rho$ and the
$\G$ action $\rho'$. The last step in the argument is to improve
the regularity of $\varphi$ which we will discuss following the
sketch of the proof of Theorem \ref{theorem:isomrigid}

We recall two definitions and another theorem from \cite{FM2}.

\begin{defn}
\label{definition:epsilonisometry} Let $\varepsilon{\geq}0$ and
$Z$ and $Y$ be metric spaces.  Then a map $h:Z{\rightarrow}Y$ is
an {\em $\varepsilon$-almost isometry} if
$$(1-{\varepsilon})d_Z(x,y){\leq}{d_Y(h(x),h(y))}{\leq}(1+{\varepsilon}){d_Z(x,y)}$$
for all $x,y{\in}Z$.
\end{defn}

\noindent The reader should note that an $\varepsilon$-almost
isometry is a bilipschitz map.  We prefer this notation and
vocabulary since it emphasizes the relationship to isometries.

\begin{defn}
\label{definition:displacement} Given a group $\Gamma$ acting on a
metric space $X$, a compact subset $S$ of $\Gamma$ and a point
$x{\in}X$.  The number $\sup_{k{\in}S}d(x,k{\cdot}x)$ is called
the $S$-displacement of $x$ and is denoted $\dk(x)$.
\end{defn}

\begin{theorem}
\label{theorem:fixedpointsimplegen} Let $\Gamma$ be a locally
compact, $\sigma$-compact group with property $(T)$ and $S$ a
compact generating set.  There exist positive constants
$\varepsilon$ and $D$, depending only on $\Gamma$ and $S$, such
that for any continuous action of $\Gamma$ on a Hilbert space
$\fh$ where $S$ acts by $\varepsilon$-almost isometries there is a
fixed point $x$; furthermore for any $y$ in $X$, the distance from
$y$ to the fixed set is not more than $D\dk(y)$.
\end{theorem}

We now sketch the proof of Theorem \ref{theorem:isomrigid} for
Theorem \ref{theorem:fixedpointsimplegen}. Given a compact
Riemannian manifold $X$, there is a canonical construction of a
Sobolev inner product on $C^k(X)$ such that the Sobolev inner
product is invariant under isometries of the Riemannian metric,
see \cite[Section 4]{FM2}. We call the completion of $C^k(X)$ with
respect to the metric induced by the Sobolev structure
$L^{2,k}(X)$. Given an isometric $\Gamma$ action $\rho$ on a
manifold $M$ there may be no non-constant $\Gamma$ invariant
functions in $L^{2,k}(X)$. However, if we pass to the diagonal
$\Gamma$ action on $X{\times}X$, then any function of the distance
to the diagonal is $\Gamma$ invariant and, if $C^k$, is in
$L^{2,k}(X{\times}X)$.

We choose a smooth function $f$ of the distance to the diagonal in
$X{\times}X$ which has a unique global minimum at $x$ on
$\{x\}{\times}X$ for each $x$, and such that any function $C^2$
close to $f$ also has a unique minimum on each $\{x\}{\times}X$.
This is guaranteed by a condition on the Hessian and the function
is obtained from $d(x,y)^2$ by renormalizing and smoothing the
function away from the diagonal. This implies $f$
 is invariant under the diagonal $\Gamma$ action defined by $\rho$.
 Let $\rho'$ be another action
$C^k$ close to $\rho$. We define a $\Gamma$ action on $X{\times}X$
by acting on the first factor by $\rho$ and on the second factor
by $\rho'$. For the resulting action $\bar \rho'$ of $\Gamma$ on
$L^{2,k}(X{\times}X)$ and every $\g{\in}S$, we show that $\bar
\rho'(k)$ is an $\varepsilon$-almost isometry and that the
$S$-displacement of $f$ is a small number $\delta$, where both
$\varepsilon$ and $\delta$ can be made arbitrarily small by
choosing $\rho'$ close enough to $\rho$. Theorem
\ref{theorem:fixedpointsimplegen} produces a $\bar \rho'$
invariant function $f'$ close to $f$ in the $L^{2,k}$ topology.
Then $f'$ is $C^{k-\dim(X)}$ close to $f$ by the Sobolev embedding
theorems and if $k-\dim(X){\geq}2$, then $f$ has a unique minimum
on each fiber $\{x\}{\times}X$ which is close to $(x,x)$. We
verify that the set of minima is a $C^{k-\dim(X)-1}$ submanifold
and, in fact, the graph of a conjugacy between the $\Gamma$
actions on $X$ defined by $\rho$ and $\rho'$.

Note that this argument yields worse regularity than we discussed
in the foliated context or than is stated in Theorem
\ref{theorem:isomrigid}.  There are considerable difficulties
involved in achieving lower loss of regularity or in producing a
$C^{\infty}$ conjugacy and we do not dwell on these here, but
refer the reader to \cite{FM2} and also to some discussion in the
next step.

\subsection*{\bf Step 4: Regularity of the conjugacy.}

We improve the regularity of $\varphi$ in three stages. First we
show it is a homeomorphism in \cite[Section 5.2]{FM3}. The key
to this argument is proving that there is a finite collection
$\g_1,{\ldots},\g_k$ of elements of $\G$ such that $\varphi$ takes
stable foliations for $\rho(\g_i)$ to stable foliations for
$\rho'(\g_i)$. If $\varphi$ were continuous, this is both easy and
classical.  In our context, we require a density point argument to
prove this, which is given in \cite[Section 5.1]{FM3}.  Once we
have that stable foliations go to stable foliations, we use this
to show that $\varphi$ is actually uniformly continuous along
central foliations and then patch together continuity along
various foliations. (In \cite[Section 5.3]{FM3}, we show how to
remove the assumption, made above, that the $\G$ action lifts from
$H/{\Lambda}$ to $H$. This cannot be done until we have produced a
continuous conjugacy.)

The next stage is to show that $\varphi$ is a finite regularity
diffeomorphism.  To show this, we show that $\varphi$ is smooth
(with estimates) along certain foliations which span the tangent
space.  This step is essentially an implementation of the method
of Katok and Spatzier described in subsection
\ref{subsection:volumepreserving}.  A few technical difficulties
occur as we need to keep careful track of estimates in the method
for use at later steps and because we need to identify ergodic
components of the measure for certain elements in the unperturbed
action.  After applying the Katok--Spatzier method, we have that
$\varphi$ is smooth along many contracting directions and smooth
along the central foliation, and then use a fairly standard
argument involving elliptic operators to show that it is actually
smooth, and even $C^{k'}$ small, for
$k'=k-1-\varepsilon-\half{\dim(H)}$.  It would be interesting
to see if one could lose less regularity at this step, for example
by a method like Journ\'{e}'s.  A key difficulty in adapting the
method of \cite{Jn} is that we only have a Sobolev estimate along
central leaves and not a uniform one.

The last stage of the argument is to show that $\varphi$ is
smooth. There are two parts to this argument.  The first is to
show that if $\rho$ and $\rho'$ are $C^k$ close, we can actually
show that $\varphi$ is $C^l$ for some $l{\geq}k$.  The main
difficulty here is obtaining better regularity in the foliated
version of Theorem \ref{theorem:isomrigid}.  This requires the use
of estimates on convexity of derivatives and estimates on
compositions of diffeomorphisms, as well as an iterative method of
constructing the semi-conjugacy $\psi$, for this we refer the
reader to \cite[Sections 6 and 7.3]{FM2}.  Once we know
we can produce a conjugacy of greater regularity, we can then
approximate $\varphi$ in the $C^l$ topology by a $C^{\infty}$ map
and smoothly conjugate $\rho'$ to a very small $C^l$ perturbation
perturbation of $\rho$. The point is to iterate this procedure
while obtaining estimates on the size of the conjugacy produced at
each step. We then show that the iteration converges to produce a
smooth conjugacy. We give here a general theorem whose proof
follows from arguments in \cite{FM2,FM3}.

It is convenient to fix right invariant metrics $d_l$ on the
connected components of $\Diff^l(X)$ with the additional property
that if $\varphi$ is in the connected component of
$\Diff^{\infty}(X)$, then $d_l(\varphi,
\Id){\leq}d_{l+1}(\varphi,\Id)$.  To fix $d_l$, it suffices to
define inner products $\<\ {,}\ \>_l$ on $\Vect^l(X)$ which satisfy
$\<V,V\>_l{\leq}\<V,V\>_{l+1}$ for $V{\in}\Vect^{\infty}(X)$. As
remarked in \cite[Section 6]{FM2}, after fixing a Riemannian
metric $g$ on $X$, it is straightforward to introduce such metrics
using the methods of \cite[Section 4]{FM2}.

\begin{defn}
\label{definition:stronglylocallyrigid} Let $\G$ be a group, $M$ a
compact manifold and assume
$\rho:\G{\rightarrow}\Diff^{\infty}(M)$.  We say $\rho$ is {\em
strongly $C^{k,l,n}$ locally rigid} if for every $\varepsilon>0$
there exists $\delta>0$ such that if $\rho'$ is an action of $\G$
on $M$ with $d_k(\rho'(g)\rho(g){\inv},\Id)<\delta$ for all
$g{\in}K$ then there exists a $C^l$ conjugacy $\varphi$ between
$\rho$ and $\rho'$ such that $d_{k-n}(\varphi,\Id)<\epsilon$.
\end{defn}

\noindent We are mainly interested in the case where $l>k$.

\begin{theorem}
\label{theorem:improvingregularity} Let $\G$ be a group, $M$ a
compact manifold and assume
$\rho:\G{\rightarrow}\Diff^{\infty}(M)$. Assume that there are
constants $n>0$ and $k_0{\geq}0$  and that for every $k>{k_0}$
there exists an $l>k$ such that $\rho$ is strongly $C^{k,l,n}$
locally rigid.  Then $\rho$ is $C^{\infty}$ locally rigid.
\end{theorem}

\noindent The proof of Theorem \ref{theorem:improvingregularity}
follows \cite[Corollary 7.2]{FM3}, though the result is not stated
in this generality there.

\subsection{KAM, implicit function theorems: work of
Damjanovich--Katok and the author.} \label{subsection:kam}

In this subsection, we discuss some new results that use more
geometric/analytic methods to approach questions of local
rigidity. These methods are entirely independent of methods using
hyperbolic dynamics and appear likely to be robustly applicable.

We begin with a theorem of Katok and Damjanovich concerning
abelian groups of toral automorphisms.  Here we consider actions
$\pi:\Za^n{\rightarrow}\Diff^{\infty}(\Ta^m)$ where $\pi(\Za^n)$
lies in $GL(m,Z)$ acting on $\Ta^m$ by linear automorphisms or
more generally where $\pi(\Za^n)$ acts affinely on $\Ta^m$. An
{\em affine factor} $\pi'$ of $\pi$ is another affine action
$\pi':\Za^n{\rightarrow}{\Ta^l}$ and there is an affine map
$\phi:\Ta^n{\rightarrow}\Ta^l$ such that
$\phi{\circ}\pi(v)=\pi'(v){\circ}\phi$ for every $v{\in}\Za^n$. We
say a factor $\pi'$ has {\em rank one} if $\pi'(\Za^n)$ has a
finite index cyclic subgroup.

\begin{theorem}
\label{theorem:damjanovichkatok}Let
$\pi:\Za^n{\rightarrow}\Aff(\Ta^m)$ have no rank one factors. Then
$\pi$ is locally rigid.
\end{theorem}

\noindent This theorem is proven by a KAM method.  One should note
that all the theorems on actions of abelian groups by toral
automorphisms stated in subsection
\ref{subsection:volumepreserving} are special cases of this
theorem.  (This is not quite literally true.  Those theorems apply
to perturbations that are only $C^1$ close, while the result
currently under discussion only applies to actions that are close
to very high order.)  It is also worthwhile to note that Theorem
\ref{theorem:damjanovichkatok} can be proven using no techniques
of hyperbolic dynamics.

We now briefly describe the KAM method.  Let $\G$ be a finitely
generated group and $\pi:\G{\rightarrow}{\Diff^{\infty}(M)}$ a
homomorphism.  To apply a KAM-type argument, define
$L:\Diff^{\infty}(M)^k{\times}\Diff^{\infty}(M){\rightarrow}\Diff^{\infty}(M)^k$
by taking
$$L(\phi_1,\ldots,\phi_k,f)=(\pi(\gamma_1){\circ}f{\circ}\phi_1{\circ}f{\inv},\ldots,\pi(\gamma_1){\circ}f{\circ}\phi_1{\circ}f{\inv})$$
where $\pi:\G{\rightarrow}\Diff^{\infty}(M)$ is a homomorphism. If
$\pi'$ is another $\G$ action on $M$ then
$L(\pi'(\gamma_1),\ldots,\pi'(\gamma_k),f)=(\Id,\ldots,\Id)$
implies $f$ is a conjugacy between $\pi$ and $\pi'$, so the
problem of finding a conjugacy is the same as finding a
diffeomorphism $f$ which solves
$L(\pi'(\gamma_1),\ldots,\pi'(\gamma_k),f)=(\Id,\ldots,\Id)$
subject to the constraint that $\pi'$ is a $\Gamma$ action.

The $KAM$ method proceeds by taking the derivative $DL$ of $L$ at
$(\pi,f)$ and solving the resulting linear equation instead
subject to a linear constraint that is the derivative of the
condition that $\pi'$ is a $\G$ action.  This produces an
``approximate solution" to the non-linear problem and one proceeds
by an iteration.  If the original perturbation is of size
$\varepsilon$ then the perturbation obtained after one step in the
iteration is of size $\varepsilon^2$ at least to whatever order
one can control the size of the solution of $DL$. This allows one
to show that the iteration converges even under conditions where
there is some dramatic ``loss of regularity", i.e. when solutions
of $DL$ are only small to much lower order than the initial data.
A standard technique used to combat this loss of regularity is to
alter the equation given by $DL$ by inserting smoothing operators.
That one can solve the linearized equation modified by smoothing
operators in place of the original linearized equation and still
expect to prove convergence of the iterative procedure depends
heavily on the quadratic convergence just described.  The main
difficulty in applying this outline is obtaining so-called tame
estimates on inverses of linearized operators.  For a definition
of a tame estimate, see following Theorem \ref{theorem:affine}.

The KAM method is often presented as a method for proving hard
implicit function theorems.  The paradigmatic theorem of this kind
is due to Hamilton \cite{Ha,Ha2}, and is used by the author in the
proof of the following theorem.  For a brief discussion of the
relation between this work and that of Katok and Damjanovich, see
the end of this subsection.

\begin{theorem}
\label{theorem:implicit} Let $\Gamma$ be a finitely presented
group, $(M,g)$ a compact Riemannian manifold and
$\pi:\G{\rightarrow}\Isom(M,g){\subset}\Diff^{\infty}(M)$ a
homomorphism. If $H^1(\G,\Vect^{\infty}(M))=0$ and
$H^2(\Gamma,\Vect^{\infty}(M))$ is Hausdorff in the tame topology,
the homomorphism $\pi$ is locally rigid as a homomorphism into
$\Diff^{\infty}(M)$.
\end{theorem}

\noindent I believe the condition on
$H^2(\Gamma,\Vect^{\infty}(M))$ should hold automatically under
the other hypotheses of the theorem.  If this is true, then one
has a new proof of Theorem \ref{theorem:isomrigid} using a result
in \cite{LZ}. There are some other infinitesimal rigidity results
that would then yield more novel applications.  For example:

\begin{theorem}
\label{theorem:irredlattices} Let $\G$ be an irreducible lattice
in a semisimple Lie group $G$ with real rank at least two. Then
for any Riemannian isometric action of $\G$ on a compact manifold
$H^1(\G,\Vect^{\infty}(M))=0$.
\end{theorem}

\noindent Theorem \ref{theorem:irredlattices} naturally applies in
greater generality, in particular to irreducible $S$-arithmetic
lattices and to irreducible lattices in products of more general
locally compact groups.

To give another result on infinitesimal rigidity, we require a
definition. For certain cocompact arithmetic lattices $\G$ in a
simple group $G$, the arithmetic structure of $\G$ comes from a
realization of $\G$ as the integer points in $G{\times}K$ where
$K$ is a compact Lie group. In this case it always true that the
projection to $G$ is a lattice and the projection to $K$ is dense.
We say a $\G$ action is {\em arithmetic isometric} if it is
defined by projecting $\G$ to $K$, letting $K$ act by $C^{\infty}$
diffeomorphisms on a compact manifold $M$ and restricting the $K$
action to $\G$.

\begin{theorem}
\label{theorem:clozel} For certain congruence lattices
$\G<SU(1,n)$, any arithmetic isometric action of $\G$ has
$H^1(\G,\Vect^{\infty}(M))=0$.
\end{theorem}

\noindent For a description of which lattices this applies to, we
refer the reader to \cite{F2}.  The theorem depends on deep
results of Clozel \cite{Cl1}.

Theorem \ref{theorem:implicit} actually follows from the
following, more general result.

\begin{theorem}
\label{theorem:affine} Let $\Gamma$ be a finitely presented group,
$M$ a compact manifold, and $\pi:\G{\rightarrow}\Diff^{\infty}(M)$
a homomorphism. If $H^1(\G,\Vect^{\infty}(M))=0$ and the sequence
$${\xymatrix{
{C^0(\G,\Vect^{\infty}(M))}\ar[r]^{d_1}&{C^1(\G,\Vect^{\infty}(M))}\ar[r]^{d_2}&{C^2(\G,\Vect^{\infty}(M))}
}}$$ \noindent admits a tame splitting then the homomorphism $\pi$
is locally rigid.  I.e. $\pi$ is locally rigid provided there
exist tame linear maps
$$V_1:C^1(\G,\Vect^{\infty}(M)){\rightarrow}C^0(\G,\Vect^{\infty}(M))$$
and
$$V_2:C^2(\G,\Vect^{\infty}(M)){\rightarrow}C^1(\G,\Vect^{\infty}(M))$$
such that $d_1{\circ}V_1+V_2{\circ}d_2$ is the identity on
$C^1(\G,\Vect^{\infty}(M))$.
\end{theorem}

\noindent Here $C^i(\G,Vect^{\infty})$ is the group of
$i$-cochains with values in $\Vect^{\infty}(M)$ and $d_i$ are the
standard coboundary maps where we have identified the cohomology
of $\G$ with the cohomology of a $K(\G,1)$ space with one vertex,
one edge for each generator in our presentation of $\G$ and one
$2$ cell for each relator in our presentation of $\G$. A map $L$
is called tame if there is an estimate of the type
$\|Lv\|_k{\leq}C_k\|v\|_{k+r}$ for a fixed choice of $r$. Here the
$\|{\cdot}\|_l$ can be taken to be the $C^l$ norm on cochains with
values in $\Vect^{\infty}(M)$. This notion clearly formalizes the
notion of being able to solve an equation with some loss of
regularity.

The proof of Theorem \ref{theorem:affine} proceeds by reducing the
question to Hamilton's implicit function theorem for short exact
sequences and is similar in outline to Weil's proof of Theorem
\ref{theorem:weil}.

Fix a finitely presented group $\G$ and a presentation of $\G$.
This is a finite collection $S$ of generators $\g_1,{\ldots},\g_k$
and finite collection $R$ of relators $w_1,{\ldots},w_r$ where
each $w_i$ is a finite word in the $\g_j$ and their inverses. More
formally, we can view each $w_i$ as a word in an alphabet on $k$
letters.  Let $\pi:\G{\rightarrow}\Diff^{\infty}(M)$ be a
homomorphism, which we can identify with a point in
$\Diff^{\infty}(M)^k$ by taking the images of the generators. We
have a complex:

\begin{equation}
\label{equation:initialsequence} \Seqinitial
\end{equation}

 \noindent Where $P$ is defined by taking
$\psi$ to $(\psi\pi(\g_1)\psi{\inv}, {\ldots},
\psi\pi(\g_k)\psi{\inv})$ and $Q$ is defined by viewing each $w_i$
as a word in $k$ letters and taking $(\psi_1,{\ldots},\psi_k)$ to
$(w_1(\psi_1,{\ldots},\psi_k),{\ldots},w_r(\psi_1,{\ldots},\psi_k))$.
To this point this is simply Weil's proof where
$\Diff^{\infty}(M)$ is replacing a finite dimensional Lie group
$H$.  Letting $\Id$ be the identity map on $M$, it follows that
$P(\Id)=\pi$ and $Q(\pi)=(\Id,{\ldots},\Id)$. Also note that
$Q{\inv}(\Id_M,{\ldots},\Id_M)$ is exactly the space of $\G$
actions.  Note that $P$ and $Q$ are $\Diff^{\infty}(M)$
equivariant where $\Diff^{\infty}(M)$ acts on itself by left
translations and on $\Diff^{\infty}(M)^k$ and
$\Diff^{\infty}(M)^r$ by conjugation.  Combining this equivariance
with Hamilton's implicit function theorem, I show that local
rigidity is equivalent to producing a tame splitting of the
sequence
\begin{equation}
\label{equation:splittingsequencepoint} \Seqpoint
\end{equation}

To complete the proof of Theorem \ref{theorem:affine}
requires that one compute $DP_{\Id}$ and $DQ_{\pi}$ in order to
relate the sequence in equation
$(\ref{equation:splittingsequencepoint})$ to the cohomology
sequence in Theorem \ref{theorem:affine}.

We remark here that the information needed to split the sequence
in Theorem \ref{theorem:affine} is quite similar to the
information one would need to apply a KAM method. This is not
particularly surprising as Hamilton's implicit function theorem is
a formalization of the KAM method.  In particular, to prove
Theorem \ref{theorem:damjanovichkatok}, one can apply Theorem
\ref{theorem:affine} using estimates and constructions from
\cite[Section 3]{DK2} to produce the required tame splitting. This
avoids the use of the explicit KAM argument in \cite[Section
4]{DK2}.

Finally, we remark that there is a theorem of Fleming that is an
analogue of Theorems \ref{theorem:implicit} and
\ref{theorem:affine} in the setting of finite, or finite Sobolev,
regularity \cite{Fl}.  This is proven using an infinite
dimensional variant of Stowe's fixed point theorems, though it has
recently been reproven by An and Neeb using a new implicit
function theorem \cite{AN}.  With either proof, this result also
has a similar condition on second cohomology.  We remark here that
due to the nature of the respective topologies on spaces of vector
fields, the condition on second cohomology in the work of Fleming
or An-Neeb is considerably stronger than what is needed in
Theorems \ref{theorem:implicit} and \ref{theorem:affine}.  As an
illustration, no version of Theorem \ref{theorem:damjanovichkatok}
can be proven using these results. This is because a cohomological
equation can have solutions with tame estimates, i.e. with some
loss of regularity, without having solutions with an estimate at
any fixed regularity.

\subsection{Further results.} \label{section:further}

In this subsection, we describe a few more recent developments
related to the results discussed so far.  These results are either
very recent or somewhat removed from the main stream of research.

The first result we discuss concerns actions of lattices in
$Sp(1,n)$ or $F_4^{-20}$ and is due to T.J.Hitchman. We state here
only a special case of his results.  For this result, we assume
that $\G$ is an arithmetic subgroup of $Sp(1,n)$ or $F_4^{-20}$ in
the standard $\Qa$ structures on those groups.  This means that
$\G$ is a finite index subgroup of the integer points in the
standard matrix representation of these groups.  This means that
in the defining representation of $Sp(1,n)$ or $F_4^{-20}$ on
$\Ra^m$ the action of $\G$ preserves the integer lattice $\Za^m$
and therefore defines an action $\rho$ of $\G$ on $\Ta^m$.

\begin{theorem}
\label{theorem:hitchman} The action $\rho$ defined in the last
paragraph is deformation rigid.
\end{theorem}

\noindent The proof proceeds in two steps.  Building a path of
$C^0$ conjugacies follows more or less as in \cite{H1}, see
subsection \ref{subsection:volumepreserving} above.  The main
novelty in \cite{Hi} is the proof that these conjugacies are in
fact smooth. Theorem \ref{theorem:hitchman} is a special case of
the results obtained in \cite{Hi}.

Another recent development should lead to a common generalization
of Theorems \ref{theorem:damjanovichkatok} and
\ref{theorem:katokspatzierRd}.  For example, one can consider
actions of an abelian subgroup $\Ra^k$ of the full diagonal group
$\Ra^{n-1}$ in $SL(n,\Ra)$ on $SL(n,\Ra)/{\Lambda}$ where
$\Lambda$ is a cocompact lattice.  In this context, Damjanovich
and Katok are developing a more geometric approach in contrast to
the analytic method of \cite{DK1,DK2}.  In \cite{DK3}, under a
natural non-degeneracy condition on the subgroup
$\Ra^k<\Ra^{n-1}$, the authors prove the cocycle rigidity result
required to generalize Theorem \ref{theorem:katokspatzierRd} to
the natural result for the $\Ra^k$ action. Here there are
``additional trivial perturbations" of the action arising from
$\Hom(\Ra^k,\Ra^{n-1})$. A rigidity theorem in this context is
work in progress, see \cite{DK3} for some discussion.

To close this section, we mention two other recent works.  The
first is a paper by Burslem and Wilkinson which investigates local
and global rigidity questions for actions of certain solvable
groups on the circle. Particularly striking is their construction
of group actions which admit $C^r$ perturbations but no $C^{r+1}$
perturbations for every integer $r$.  The second is a paper by
M.Einsiedler and T.Fisher in which the method of proof of Theorem
\ref{theorem:ksabelian} is extended to affine actions of $\Za^d$
where the matrices generating the group action have non-trivial
Jordan blocks.  For perturbations of the group action which are
close to very high order, this result follows from Theorem
\ref{theorem:damjanovichkatok}, but in \cite{EF} the result only
requires that the perturbation be $C^1$ close to the original
action.

\section{Directions for future research and conjectures}
\label{section:future}

In this section, I mention a few conjectures and point a few
directions for future research.  These are particularly informed
by my taste.

\subsection{Actions of groups with property $(T)$.}

Lattices in $SP(1,n)$ and $F_4^{-20}$ share many of the rigidity
properties of higher rank lattices.  In light of Theorems
\ref{theorem:main}, \ref{theorem:isomrigid} and
\ref{theorem:hitchman}, it seems natural to conjecture:

\begin{conjecture}
\label{conjecture:TJandI} Let $G$ be a semi-simple Lie group with
no compact factors and no simple factors isomorphic to $SO(1,n)$
or $SU(1,n)$, and let $\G<G$ be a lattice.  Then any volume
preserving generalized quasi-affine action of $G$ or $\G$ on a
compact manifold is locally rigid.
\end{conjecture}

\noindent Note that for $G$ with no rank one factors and
quasi-affine actions, this is just Theorem \ref{theorem:main}. To
apply the methods of \cite{FM1,FM2,FM3} in the setting of
Conjecture \ref{conjecture:TJandI} there are essentially three
difficulties:

\begin{enumerate}
\item If the action is generalized quasi-affine and not
quasi-affine, then one cannot use the construction of the
Margulis--Qian cocycle described above.  This is easiest to see for
a generalized affine action $\rho$ on some
$K{\backslash}H/{\Lambda}$.  The action $\rho$ lifts to
$H/{\Lambda}$, but a perturbation $\rho'$ need not.  If $K$ is
finite, this difficulty can be overcome by passing to a subgroup
of finite index $\G'<\G$ and arguments in \cite{FM3} can be used
prove rigidity of $\G$ from rigidity of $\G'$.  If $K$ is compact
and connected, this is a genuine and surprisingly intractable
difficulty.

\item Proving a version of Zimmer's cocycle super-rigidity for
groups as in the assumptions of the conjecture.  Partial results
in this direction were obtained by Corlette--Zimmer and
Korevaar--Schoen, but their results all require hypotheses that are
obviously restrictive or simply difficult to verify.  Very
recently, the author and T.J.Hitchman have proven a complete
version of cocycle super-rigidity, at least for so-called
$L^2$-cocycles.  In light of our work, this difficulty is already
overcome.

\item  Replace the method of Katok--Spatzier in the proof that the
conjugacy is actually smooth.  The proof of Theorem
\ref{theorem:hitchman} gives some progress in this direction, but
there are significant technical difficulties to overcome in
applying Hitchman's methods at this level of generality.  There is
some progress on this question by Gorodnik, Hitchman and Spatzier.
\end{enumerate}

\noindent The author and T.J.Hitchman have another approach to
Conjecture \ref{conjecture:TJandI} based on Theorem
\ref{theorem:affine}, some estimates proven in \cite{FH1}, and
using heat flow and those estimates to produce a tame splitting of
the short exact sequence in Theorem \ref{theorem:affine}. It is
not yet clear how generally applicable this method will be.

The following question is also interesting in this context:

\begin{question}
\label{question:Tboundaries} Let $G$ and $\G$ be as in Conjecture
\ref{conjecture:TJandI}.  Let $\rho$ be a non-volume preserving
affine action of $\G$ on a compact manifold $M$.  When is $\rho$
locally rigid?
\end{question}

\noindent By the work of Stuck discussed at the end of Section
\ref{subsection:boundaries}, it is clear that local rigidity will
not hold in full generality here.  In particular for the product
of the action of $\G$ on the boundary $G/P$ with the trivial
action of $\G$ on any manifold, there is already strong evidence
against local rigidity. We give two particularly interesting
special cases where we expect local rigidity to occur. The first
is to take a lattice $\G$ in $SP(1,n)$ and ask if the $\G$ action
on the boundary $SP(1,n)/P$ is locally rigid. This is analogous to
Theorems \ref{theorem:ksboundary} and \ref{theorem:kanai2}.
Another direction worth pursuing is to see if the action of say
$Sl(n,\Za)$ on ${\Ra}P^{n+k}$ is locally rigid for $n{\geq}3$ and
any $k{\geq}0$. For cocompact lattices instead of $SL(n,\Za)$ and
$k=-1$ it is Theorem \ref{theorem:ksboundary}. One can ask a wide
variety of similar questions for both compact and non-cocompact
lattices acting on homogeneous spaces that are ``larger than" any
natural boundary for the group as long as one avoids settings in
which Stuck's examples can occur.   We remark that in some
instances, partial results for analytic perturbations can be
obtained by Zeghib's method \cite{Zg}.

We remark here that the action of $SL(n, \mathbb Z)$ on $\mathbb
T^n$ is not locally rigid in $\Homeo(\mathbb T^n)$ by a
construction of Weinberger.  It would be interesting to understand
local rigidity in low regularity for other actions.

\subsection{Actions of irreducible lattices in products.}
\label{subsection:irredlattices}

In this subsection, we formulate a general conjecture concerning
local rigidity of actions of irreducible lattices in products. We
begin by making a few remarks on other rigidity properties of
irreducible lattices and by describing a few examples where
rigidity might hold, as well as some examples where it does not.

Rigidity properties of irreducible lattices have traditionally
been studied together with rigidity of lattices in simple groups,
and irreducible lattices enjoy many of the same rigidity
properties.  We list a few here to motivate our conjectures on
rigidity of actions of these lattices.  The properties we list
also rule out certain trivial constructions of perturbations and
deformations of actions.  In the following, $\G$ is an irreducible
lattice in $G=(\prod_{i=1}^kG_i)/Z$ where each $G_i$ is a
noncompact semi-simple Lie group and $Z$ is a subgroup of the
center of $\prod_{i=1}^kG_i$. Many of these results hold more
generally, see below.

\smallskip
{\bf Properties of irreducible lattices:}
\begin{enumerate}
\item There are no non-trivial homomorphisms $\G{\rightarrow}\Za$
(and therefore no non-trivial homomorphisms to any abelian or
non-abelian free group).

\item All linear representations of $\G$ are classified.  In
particular, given any representation $\rho$ of $\G$ into $GL(V)$,
where $V$ is a finite dimensional vector space, then
$H^1(\G,V)=0$.

\item All normal subgroups of $\G$ are either finite or finite
index.
\end{enumerate}
\noindent The first property, for $\G$ cocompact, was originally
proven by Bernstein and Kazhdan \cite{BK}. Some special cases of
this result were proven earlier by Matsushima and Shimura
\cite{MS}. All other properties are originally due to Margulis,
see \cite{Ma3} and historical references there. The original
proofs all go through with only minor adaptations if some of the
$G_i$ are replaced with the $k$-points of a $k$-algebraic group
over some other local field $k$.  These properties have been shown
to hold for appropriate classes of lattices in even more general
products of locally compact groups by Bader, Monod and Shalom,
\cite{Sm,Md1,Md2,BSh}.

\smallbreak
\noindent{\bf Example $1$:} Let $\G=SL(2,\Za[\sqrt{2}])$.  We
embed $\G$ in $SL(2,\Ra){\times}SL(2,\Ra)$ by taking $\g$ to
$(\g,\sigma(\g))$ where $\sigma$ is the Galois automorphism of
$\Qa(\sqrt{2})$ taking $\sqrt{2}$ to $-\sqrt{2}$.  This embedding
defines an action of $\G$ on $\Ta^4$ where $\Ta^4=\Ra^4/\Za^4$
where we identity $\Za^4$ with the image in $\Ra^4$ of
$\Za[\sqrt{2}]^2$ via the embedding $v{\rightarrow}(v,\sigma(v))$.
We first note that the list of properties given above imply one
can prove deformation rigidity of this action using Hurder's
argument from \cite{H1} to produce a continuous conjugacy and
using the method of Katok--Spatzier \cite{KS1,KS2} to show that the
conjugacy is smooth.  Anatole Katok has suggested one might be
able to show local rigidity of this action by using the methods in
\cite{KL1}. This example is just the first in a large class of
Anosov actions of irreducible lattices, all of which should be
locally rigid.  We leave the general construction to the
interested reader.

\smallbreak
\noindent{\bf Example $2$:} We take the action of
$\G=SL(2,\Za[\sqrt{2}]$ and let $\G$ act on
$\Ta^5=\Ta^4{\times}\Ta^1$ by a diagonal action where the action
on $\Ta^4$ is as defined in Example $1$ and the action on $\Ta^1$
is trivial.  Once again, this is merely the first example of a
large class of partially hyperbolic actions of irreducible
lattices.  In this instance, the central foliation for the $\G$
action consists of compact tori.  For this type of example, many
of the argument of \cite{NT1,NT2,T} carry over, but the fact that
$\G$ does not have property $(T)$ prevents one from using those
outlines to prove local rigidity.  On the other hand, the methods
of \cite{NT1} can be adapted to prove deformation rigidity again
replacing their argument for smoothness of the conjugacy by the
Katok--Spatzier method. We remark that it is also possible to give
many examples of actions of irreducible lattices in products where
the central foliation is by dense leaves, and the methods of
Nitica and Torok cannot be applied. See discussion below for the
obstructions to applying the methods of \cite{FM1,FM2,FM3}.

\smallbreak
\noindent{\bf Example $3$:} We now give an example of a family of
actions which extend to an action of $G=G_1{\times}G_2$.  Let $H$
be a simple Lie group with $G<H$, for example, $H=SL(4,\Ra)$ or
$H=Sp(4,\Ra)$ and $\Lambda<H$ a cocompact lattice.  Then both $G$
and $\G$ act by left translations on $H/{\Lambda}$.

\smallbreak
\noindent{\bf Example $4$:}We end with a family of examples for
which there exists a large family of deformations.  Let $\G$ be as
above and let $\Lambda<SL(2,\Ra)=G_1$ be an irreducible lattice.
Then $\G$ acts by left translations on $M=SL(2,\Ra)/\Lambda$. Call
this action $\bar \rho$. I do not currently know whether it is
possible to deform this action, but one can use this action to
build actions with perturbations on a slightly larger manifold.
Let $\G$ act trivially on any manifold $N$ and take the diagonal
action $\rho$ on $M{\times}N$. It is well-known that there exists
a non-trivial homomorphism $\sigma:\Lambda{\rightarrow}\Za$. There
is also a standard construction of a cocycle
$\alpha:G_1{\times}{G_1/{\Lambda}{\rightarrow}\Lambda}$ over the
$G_1$ action on $G_1/{\Lambda}$.  The cocycle is defined by taking
a fundamental domain $X$ for $\Lambda$ in $G_1$, identifying
$G_1/\Lambda$ with $X$ and letting $\alpha(g,x)$ be the element of
$\Lambda$ such that $gx\alpha(g,x){\inv}$, as an element of $G_1$,
is in $X$.  Taking any vector field $V$ on $N$ and let
$s_{\varepsilon}$ be the $\Za$ action on $N$ defined by having $1$
act by flowing to time $\varepsilon$ along $V$.  We then define a
$\G$ action on $M{\times}N$ by taking
$$\rho_{\varepsilon}(\g)(m,n)=(\bar
\rho(\g)m,s_{\varepsilon}(\sigma(\alpha(g,m))(n)).$$ \noindent We
leave it to the interested reader to show that the actions
$\rho_{\varepsilon}$ are not conjugate to $\rho$ for essentially
any choice of $V$.
\medbreak

The key point in example $4$, which is not present in the first
three examples, is that the action factors through a projection of
$\G$ into a simple factor of $G$.  Motivated by the examples so
far, by the results in Theorem \ref{theorem:irredlattices}, and by
analogy with results on actions of higher rank abelian groups, we
make the following definition and conjecture:

\begin{defn}
\label{definition:rankone} Let $G=G_1{\times}{\cdots}{\times}G_k$
be a semisimple Lie group where all the $G_i$ are non-compact. Let
$\G<G$ be an irreducible lattice.  Let $\rho$ be an affine action
of $\G$ on some $H/{\Lambda}$ where $H$ is a Lie group and
${\Lambda}<H$ is a cocompact lattice.  Then we say $\rho$ has {\em
rank one factors} if there exists \begin{enumerate} \item an
action $\bar \rho$ of $\G$ on a some space $X$ which is a factor
of $\rho$ \item and a rank one factor $G_i$ of $G$
\end{enumerate}
such that $\bar \rho$ is the restriction of a $G_i$ action.  I.e.
$\G$ acts on $X$ by projecting $\G$ to $G_i$ and restricting a
$G_i$ action.
\end{defn}

\begin{conjecture}
\label{conjecture:irredlattices} Let $G,\G,H,\Lambda$ and $\rho$
be as in Definition \ref{definition:rankone}.  Then if $\rho$ has
no rank one factors $\rho$ is locally rigid.
\end{conjecture}

\begin{remark}
It is a consequence of Ratner's measure
rigidity theorem, see \cite{R,Sh,W}, that any rank one factor of
an affine action for these groups is in fact affine. This implies
that any rank one factor is a left translation action on some
$H'/{\Lambda'}$.  So a special case of the conjecture is that any
affine $\G$ action on a torus or nilmanifold is locally rigid.
\end{remark}

There is a variant of Conjecture \ref{conjecture:irredlattices}
for $G$ actions.  We recall that a measure preserving action of
$G=G_1{\times}{\cdots}{\times}G_k$ action is {\em irreducible} if
each $G_i$ acts ergodically.  We extend this notion to non-ergodic
$G$ actions by saying that the action is {\em weakly irreducible}
if every ergodic component of the volume measure for the action of
any $G_i$ is an ergodic component of the volume measure for the
action of $G$.

\begin{conjecture}
\label{conjecture:irredactions} Let $G,H,\Lambda$ be as in
Definition \ref{definition:rankone} and let $G<H$, so that we have
a left translation action $\rho$ of $G$ on $H/\Lambda$.  Then if
$\rho$ is weakly irreducible, $\rho$ is locally rigid.
\end{conjecture}

The relation of the two conjectures follows from the following
easy lemma.

\begin{lemma}
\label{lemma:induction} Let $G,\G,H,\Lambda$ and $\rho$ be as in
Definition \ref{definition:rankone}, then the $\G$ action on
$H/{\Lambda}$ has no rank one factors if and only if the induced
$G$ action on $(G{\times}H/{\Lambda})/\G$ is weakly irreducible.
\end{lemma}

To prove the lemma requires both some algebraic untangling of
induced actions and a use of Margulis' super-rigidity theorem to
describe affine actions of $G$ and $\G$ along the lines of
\cite[Theorems 6.4 and 6.5]{FM1}. We leave this as an exercise for
the interested reader. It is fairly easy to check the lemma for
any particular affine $\G$ action.

We end this subsection by pointing out the difficulty in
approaching this conjecture by means of the methods of
\cite{FM1,FM2,FM3}.  A central difficulty is that the lattices in
question do not have property $(T)$ and so the foliated
generalization of Theorem \ref{theorem:isomrigid} does not apply.
However, even in the case of weakly hyperbolic actions, there are
significant difficulties.  To begin the argument, one would like
to apply cocycle superrigidity to the Margulis--Qian cocycle. To do
this requires the existence of an invariant measure which is
usually established using property $(T)$ by an argument of Seydoux
\cite{Sy}.  In this setting, where property $(T)$ does not hold,
one might try instead to use the work of Nevo and Zimmer,
\cite{NZ1,NZ2,NZ3}, but there are non-trivial difficulties here as
well.  One cannot apply their theorems without first showing that
the perturbed action satisfies some irreducibility assumption.
Even if one were to obtain an invariant measure, the precise form
of cocycle superrigidity required is not known for products of
rank one groups or their irreducible lattices. And the strongest
possible forms of cocycle superrigidity in this context again
require a kind of irreducibility of the perturbed action. So to
proceed by this method one would need to show that perturbations
of the actions in Conjecture \ref{conjecture:irredlattices} and
\ref{conjecture:irredactions} still satisfied some irreducibility
conditions.  This seems quite difficult.  It may also be possible
to approach these questions by using Theorem \ref{theorem:affine},
but even proving that the relevant cohomology groups vanish seems
subtle.

\subsection{Other questions and conjectures.}

We end this article by discussing some other questions and
conjectures.

In the context of Theorem \ref{theorem:clozel} it is interesting
to ask if isometric actions of lattices in $SU(1,n)$ are locally
rigid.  For some choices of lattice, the answer is trivially no.
Namely some cocompact lattices in $SU(1,n)$ have homomorphisms
$\rho$ to $\Za$ \cite{Ka,BW}, and so have arithmetic actions with
deformations provided the centralizer $Z$ of $K$ in
$\Diff^{\infty}(M)$ is non-trivial. Having centralizer allows one
to deform the action along the image of the homomorphism
$\rho{\circ}\sigma_t:F{\rightarrow}Z$ where
$\sigma_t:\Za{\rightarrow}Z$ is any one parameter family of
homomorphisms.  It seems reasonable to conjecture:

\begin{conjecture}
\label{conjecture:su(1,n)} Let $\rho$ be an arithmetic isometric
action of a lattice in $SU(1,n)$.  Then any sufficiently small
perturbation of $\rho$ is of the form described in the previous
paragraph.
\end{conjecture}

This conjecture is in a certain sense an infinite dimensional analogue
of work of Goldman--Millson and Corlette \cite{Co1,GM}. Another
conjecture concerning complex hyperbolic lattices, for which work
of Yue provides significant evidence \cite{Yu}, is:

\begin{conjecture}
\label{conjecture:boundarysu(1,n)} Is the action of any lattice in
$SU(1,n)$ on the boundary of complex hyperbolic space locally
rigid?
\end{conjecture}

There are also many interesting questions concerning the failure
of local rigidity for lattices in $SO(1,n)$.  The only rigidity
theorem we know of in this context is Kanai's, Theorem
\ref{theorem:kanai1}, and it would be interesting to extend
Kanai's theorem to non-uniform lattices.  In \cite{F1,F3} various
deformations of lattices in $SO(1,n)$ are constructing for affine
and isometric actions.  These constructions both adapt the bending
construction of Johnson and Millson, \cite{JM}.  It seems likely
that in some cases one should be able to prove results concerning
the structure of the representation space and, in particular, to
show that it is ``singular" in an appropriate sense. See \cite{F3}
for more discussion.

Two other paradigmatic examples of large groups are the outer
automorphism group of the free group, $\Out(F_n)$, and the mapping
class group of a surface $S$, $MCG(S)$.  These groups do not admit
many natural actions on compact manifolds, but there are some
natural interesting actions quite analogous to those we have
already discussed. For $MCG(S)$, the question we raise here is
already raised in \cite{La}.  The actions we consider are
``non-linear" analogues of the standard actions of $SL(n,\Za)$ on
$\Ta^n$ and $SP(2n,\Za)$ on $\Ta^{2n}$.  The spaces acted upon are
moduli spaces of representations of either the free group or the
fundamental group of a surface $S$, where the representations take
values in compact groups.  More precisely, we have an action of
$\Out(F_n)$ on $\Hom(F_n,K)/K$ and an action of $MCG(S)$ on
$\Hom(\pi(S),K)/K$ where $K$ is a compact group.  It is natural to
ask whether these actions are locally rigid, though the meaning of
the question is somewhat obscured by the fact that the
representation varieties are not smooth.  For $K=S^1$, one obtains
actions on manifolds, and in fact tori, and one might begin by
considering that case.

We end with a question motivated by the recent work of Damjanovic
and Katok.  We only give a special case here. Let $G$ be a real
split, simple Lie group of real rank at least two.  Let
$\G<G{\times}G$ be an irreducible lattice.  Let $K<G$ be a maximal
compact subgroup and view $K$ as a subgroup of $G{\times}G$ by
viewing it as a subgroup of the second factor.  The quotient
$K{\backslash}(G{\times}G)/{\Gamma}$ has a natural $G$ action on
the left on the first factor.  We can restrict this action to the
action of a maximal split torus $A$ in $G$.  Note that $A$ is
isomorphic to $\Ra^d$ for some $d{\geq}2$.

\begin{question}
\label{question:rsplit} Is the action of $\Ra^d$ described in the
paragraph above locally rigid?
\end{question}

\end{document}